\documentstyle[fleqn,11pt]{article}

\catcode`\@=11

 \def\AMSTeXfeatures{\Plainheads 
   \let\current@vert=\AMS@vert}

 \def\Plainheads{\sh@ftdiam=0.05em
   \getlabeldims
   \let\vshaftfill=\plnvsolidfill
   \let\hshaftfill=\plnhsolidfill
   \let\th@rhead=\plnrhead
   \let\th@lhead=\plnlhead
   \let\th@dnhead=\plndnhead
   \let\th@uphead=\plnuphead}
 
 \def\glet{\global\let}

 \def\LaTeXfeatures{\catcode`\@=11
   \ifx\@clnwd\undefined \nol@g
      \input ltxcode.tex \dol@g \fi
   \ltxheads \let\current@vert=\new@vert
   \providelto \catcode`\@=\active}

 \def\nol@g{\def\wlog{\edef\garbage}}
 \def\dol@g{\let\wlog=\wl@g} \let\wl@g=\wlog
 \nol@g 

 \newbox\ltobox
 \def\providelto{{\setbox\z@=
   \hbox{$\to$}\minharrlen=\wd\z@
   \global\setbox\ltobox=\hbox{$\activeat>>>$}}
   \def\lto{\mathrel{\copy\ltobox}}}

 \def\ltxheads{\sh@ftdiam=\@wholewidth
   \getlabeldims
   \let\vshaftfill= \ltxvsolidfill
   \let\hshaftfill=\ltxhsolidfill
   \let\th@rhead=\ltxrhead
   \let\th@lhead=\ltxlhead
   \let\th@dnhead=\ltxdnhead
   \let\th@uphead=\ltxuphead}
 {\catcode`\@=\active
   \gdef@#1{\csname #1\string@at\endcsname}
   \glet\activeat=@}
 \def\def@#1{\expandafter\def\csname #1@at\endcsname}

 \def@>#1>#2>{\@rrow R{#1}{#2}}
 \def@<#1<#2<{\@rrow L{#1}{#2}}
 \def@ V#1V#2V{\@rrow V{#1}{#2}}
 \def@ A#1A#2A{\@rrow A{#1}{#2}}
 \def@/#1/#2/#3/{\@rrow{#1}{#2}{#3}}
 \def@.{\ifodd\row\ifmmode\noharrow
     \else\leavevmode.\spacefactor3000 \fi
   \else\novarrow\fi}
 \def@={\ifodd\row\harrow\hequalfill{}{}%
   \else\varrow\vequalfill{}{}\fi}
 \def@:#1{\ifx=#1\harrow\deffill{}{}%
   \else\leavevmode\null:#1\fi}
 \def@|{\current@vert}
  \def\AMS@vert{\varrow\vequalfill{}{}}
  \def\new@vert#1|#2|{\ifodd\row
   \let\nextarrow\vertexvarrow
   \else\let\nextarrow\varrow\fi
   \nextarrow\vshaftfill{#1}{#2}}
 \def@-{\ifmmode\let\next\hl@ne
   \else\let\next\AMSatdash \fi \next}
  \def\hl@ne#1-#2-{\harrow\hshaftfill{#1}{#2}}
  \def\AMSatdash{\let\next\relax\leavevmode
    \def\next@{\ifx\next-%
      \def\next-{\futurelet\next\nextii@}%
     \else\def\next{\hbox{-}}\fi\next}%
    \def\nextii@{\ifx\next-\def\next-{\hbox{---}}%
      \else\def\next{\hbox{--}}\fi\next}%
    \futurelet\next\next@}
 \def@(#1){\tweenarrows{#1}}
 \def@[#1]{\setsp@n#1\relax\activeat}
 \def\fiberbox{\hbox{$\vcenter{\hr@le\hbox{\vr@le
   \kern1ex\vbox{\kern1.2ex}\vr@le}\hr@le}$}}
  \def\hr@le{\hrule height \sh@ftdiam}
  \def\vr@le{\vrule width \sh@ftdiam}
 \def@+#1+#2+#3+{\ifodd\row \harrow{#1}{#2}{#3}%
   \else \varrow{#1}{#2}{#3}\fi}

 \def\Dnarrfill{\vequalfill\Dnhe@d}
 \def\Uparrfill{\Uphe@d\vequalfill}
 
 \def\ontofill{\rtarrfill\kern-0.3em 
   \th@rhead\kern 0.3em} 

 \def\rtarrfill{\hshaftfill\th@rhead}
 \def\ltarrfill{\th@lhead\hshaftfill}
 \def\dnarrfill{\vshaftfill\th@dnhead}
 \def\uparrfill{\th@uphead\vshaftfill}
 \def\hequalfill{\plnhfill=}
 \def\deffill{:\plnhfill=}
 \def\plnvextfill#1{\setbox\z@
   \hbox{\the\textfont3 #1}%
   \dimen@=\dp\z@\advance\dimen@\ht\z@
   \copy\z@ \kern-\dimen@ 
   \cleaders\copy\z@ \vfill
   \kern-\dimen@ 
   \box\z@}
 \def\plnhfill#1{$\m@th\mkern-1.5mu\mathord#1\mkern-6mu
    \cleaders\hbox{$\mkern-2mu\mathord#1\mkern-2mu$}\hfill
    \mkern-6mu\mathord#1\mkern-1.5mu$}
 \def\vequalfill{\plnvextfill{\char'167}}
 \def\plnvsolidfill{\plnvextfill{\char'077}}
 \def\plnhsolidfill{\plnhfill-}
 \def\ltxhsolidfill{\leaders\hrule height\topofshaft depth\botofshaft
   \hfill}
 \def\ltxvsolidfill{\leaders\vrule width\sh@ftdiam\vfill}
 \def\hdashfill{\hd@sh\wd@sh
   \xleaders \hbox{\wd@sh\hd@sh\wd@sh}\hfill
   \wd@sh\hd@sh}
 \def\vdashfill{\vd@sh\wd@sh
   \xleaders \vbox{\wd@sh\vd@sh\wd@sh}\vfill
   \wd@sh\vd@sh}
 \def\dashed{\ifinmeasureCD\else
    \ifodd\row\option{\let\hshaftfill=\hdashfill}%
   \else\option{\let\vshaftfill=\vdashfill}\fi\fi}

 \newdimen\CDstrutht  \newdimen\CDstrutdp
 \newdimen\CDstrutlen \CDstrutlen=\CDstrutht
   \advance\CDstrutlen by \CDstrutdp

 \def\CDstrut{\vrule
   height \ifnum\row=1 \z@\else\CDstrutht \fi
   depth \ifnum\row=\numrows \z@ \else\CDstrutdp \fi
   width\z@}

 \newdimen\CDarrsurr \CDarrsurr=0.375em
 \newdimen\CDdashlen
 \newdimen\CDvarrlen \CDvarrlen=1.5\baselineskip
 \newdimen\minharrlen 
  \setbox\z@\hbox{$\longrightarrow$} \minharrlen=\wd\z@
 \newdimen\minCDharrlen \minCDharrlen=2.5em 
\newdimen \minc@lwd
\def\findminc@lwd{\minc@lwd=2\CDarrsurr
  \advance\minc@lwd\minCDharrlen}

 \newdimen\sh@ftdiam

 \newdimen\labelsurr \labelsurr=1.25 em

\newcount\sp@ncnt \sp@ncnt=\@ne
\newcount\sp@ncnt@ \sp@ncnt@=\@ne
\newdimen\@rrwd \newdimen\@rrdp

 \def\adjustbot#1{\option{\advance\@rrdp#1\relax}}
 
\def\pushvertex#1{\global\p@shlen#1\relax
   \global\let\maybepush=\dopush}

 \newdimen\p@shlen \p@shlen=\z@

 \let\maybepush=\relax
 \def\dopush{\ifinmeasureCD 
   \advance\locdimen by -\p@shlen 
   \else\advance \@rrwd by -\p@shlen \fi 
   \global\let\maybepush=\relax \global\p@shlen=\z@\relax}

 \def\span@ne{\global\sp@ncnt=\@ne\relax}
 \def\setsp@n#1#2{\global\sp@ncnt=#1\relax
   \ifx\relax#2\relax\else\global\sp@ncnt@=#2\relax\fi}

 \def\plnrhead{\llap{$\rightarrow\mkern-1.5mu$}}
 \def\plnlhead{\rlap{$\mkern-1.5mu\leftarrow$}}

 \def\clap#1{\hbox to \z@{\hss #1\hss}}

 \def\plndnhead{\hbox{\the\textfont3 \char'171}}
 \def\plnuphead{\hbox{\the\textfont3 \char'170}}
 \def\Dnhe@d{\hbox{\the\textfont3 \char'177}}
 \def\Uphe@d{\hbox{\the\textfont3 \char'176}}

 \def\ltxrhead{\raise\@xisheight
   \llap{\smash{\@linefnt\@getrarrow(1,0)}}}
 \def\ltxlhead{\raise\@xisheight
   \rlap{\@linefnt\@getlarrow(-1,0)}}
 \def\ltxuphead{\setbox\z@=\rlap{%
   \kern\@halfwidth\@linefnt\char'66}%
   \copy\z@\kern-\ht\z@}
 \def\ltxdnhead{\setbox\z@=\rlap{%
   \kern\@halfwidth\@linefnt\char'77}%
   \ht\z@=\z@\box\z@}

 \def\wd@sh{\kern0.5\CDdashlen}
 \def\hd@sh{\vrule height\topofshaft depth\botofshaft
    width\CDdashlen}
 \def\vd@sh{\hrule height\CDdashlen
   depth\z@ width\sh@ftdiam}

\def\xylist{14{3434}13{2414}12{1723}%
  23{1413}34{1153}11{0867}43{0707}%
  32{0580}21{0414}31{0291}41{0}}
\newcount\tgtcnt@
\def\find@xyargs{\dimen@=\@rrdp
  \advance\dimen@ by \CDstrutlen
  \tgtcnt@=\dimen@ \dimen@=\@rrwd 
  \divide\dimen@ by \@m 
  \divide \tgtcnt@ by \dimen@ 
  \expandafter\testxy\xylist\relax
  \unitlength=\@xarg\@rrdp
  \divide\unitlength by\@yarg\relax}
\def\testxy#1#2#3{\ifnum\tgtcnt@>#3
    \@xarg=#1\relax \@yarg=#2\relax
    \let\next=\ignorerest
  \else\let\next\testxy\fi\next}
\def\ignorerest#1\relax{\relax}

\let\scalefactor=\@ne
\def\SWarrow{\find@xyargs\vector
  (-\@xarg,-\@yarg)\scalefactor\hskip-\wd\@linechar}
\def\NWarrow{\find@xyargs\vector
  (-\@xarg,\@yarg)\scalefactor\hskip-\wd\@linechar}
\def\NEarrow{\find@xyargs\vector
  (\@xarg,\@yarg)\scalefactor}
\def\SEarrow{\find@xyargs\vector
  (\@xarg,-\@yarg)\scalefactor}
\def\rightupline{\find@xyargs\@linelen=\scalefactor
     \unitlength\@sline}
\def\rightdownline{\find@xyargs\@yarg=-\@yarg\relax
     \@linelen=\scalefactor\unitlength\@sline}

\def\Sim{\ifodd\row\setbox\z@=\hbox{$\sim$}\dimen@=\ht\z@
 \advance\dimen@ by -\@xisheight
  \vbox{\box\z@\kern-\@xisheight\kern\dimen@}%
  \else\hbox{$\wr$}\fi}

\def\harrow#1#2#3{\inmeasureCDtrue\findminarrwd
  {#2}{#3}{\sp@ncnt\minharrlen}\inmeasureCDfalse\span@ne
  \mathrel{\hbox{\options\hplace{#1}\ulabel{#2}\dlabel{#3}}}}

\def\noharrow{\harrow\hfill{}{}}
\def\vertexvarrow#1#2#3{\findarrdp \@rrwd=\z@ \setsp@n\@ne\@ne
  \vbox to \z@{\kern-1.2\CDstrutht
  \rlap{\options\vplace{#1}\llabel{#2}\rlabel{#3}}\vss}}

\newif\ifinmeasureCD
\def\measurelabel#1{\setbox\z@
  \hbox{$\scriptstyle#1\kern\labelsurr$}%
  \ifdim\wd\z@>\@rrwd \@rrwd=\wd\z@\fi}
\def\findminarrwd#1#2#3{\@rrwd=#3\relax
   \measurelabel{#1}\measurelabel{#2}}
\def\findCDarrwd#1#2{\@rrwd=\minCDharrlen
   \measurelabel{#1}\measurelabel{#2}%
  }

\newcount\row \row=\@ne \newcount\col \col=\@ne 
 \newcount\numrows 
\numrows=\@ne
 \newcount\numcols
\newcount\arrspan \newdimen\vrtxhalfwd  \newbox\tempbox

\def\DANABUG{\advance\col by \@ne
 \@rrwd=\minCDharrlen
  \advance\@rrwd by \vrtxhalfwd
  \advance\@rrwd by \CDarrsurr
  \ifnum\col>\numcols \numcols=\col
     \newlocdimen{col\the\col}\locdimen=\@rrwd 
  \else \ifdim\@rrwd>\c@l \c@l=\@rrwd\fi\fi}

\def\drop#1\\{
  \findvrtxhalfsum\DANABUG\advance\row by 2 \measureinit}

\def\measureinit{\col=\@ne \vrtxhalfwd=-\CDarrsurr\arrspan=\@ne\@rrwd=\z@
   \setbox\tempbox=\hbox\bgroup$}
\def\measure{
  \let\harrow\measureCDarrow
  \let\CDCR=\measureCR 
   \findminc@lwd 
  \inmeasureCDtrue
  \row=\@ne \numcols=\z@ \measureinit}

\def\endmeasure{\findvrtxhalfsum\DANABUG
  \numrows=\row 
  \inmeasureCDfalse}

\def\newlocdimen#1{\advance\dimenc@unt by \@ne
  \ifnum\dimenc@unt<\insc@unt
     \else\errmessage{No room for the CD}\fi
  \dimendef\locdimen=\dimenc@unt
  \expandafter\dimendef\csname#1\endcsname=\dimenc@unt}

 \def\r@wc@l{\csname row\the\row col\the\col\endcsname}
 \def\c@l{\csname col\the\col\endcsname}

 \def\findvrtxhalfsum{$\egroup
  \newlocdimen{row\the\row col\the\col}
  \locdimen=\vrtxhalfwd 
  \vrtxhalfwd=0.5\wd\tempbox 
  \advance\vrtxhalfwd by \CDarrsurr
  \advance\locdimen by \vrtxhalfwd 
  \advance\@rrwd by \locdimen 
  \maybepush
  \divide\@rrwd by \arrspan\relax
  \ifdim\@rrwd<\minc@lwd
    \ifnum\col>\@ne \@rrwd=\minc@lwd\fi \fi
  \loop 
    \ifnum\col>\numcols \numcols=\col
       \newlocdimen{col\the\col}
       \locdimen=\@rrwd 
    \else \ifdim\@rrwd>\c@l \c@l=\@rrwd\fi \fi
   \ifnum\arrspan>\@ne
      \advance\arrspan by -1 \advance\col by \@ne
  \repeat }

 \def\measureCDarrow#1#2#3{\findvrtxhalfsum
   \arrspan=\sp@ncnt\relax\global\sp@ncnt=1\relax
   \advance\col by \@ne
   \findCDarrwd{#2}{#3}%
   \setbox\tempbox=\hbox\bgroup$}

 \newcount\dr@tn \dr@tn=\z@
 \def\locate#1:#2{\ifinmeasureCD\else
   \count@=-#1
   \multiply\count@ by 2
   \advance\count@ by #2
   \dimen@=\count@\@rrwd
   \ifnum\dr@tn=\@ne\relax \else\dimen@=-\dimen@ \fi
   \dimen@i=\@rrdp
   \ifnum\dr@tn>\z@\advance\dimen@i by \CDstrutlen \fi
   \dimen@i=\count@\dimen@i
   \count@=#2 \multiply\count@ by 2
   \divide\dimen@ by \count@
   \divide\dimen@i by \count@
   \lift\dimen@i\nudge\dimen@\fi}

\def\betweenCDrows{\advance\row by \@ne \col=\@ne
\options}

\def\hbegin{\hbox\bgroup\kern\c@l \kern-\r@wc@l$}
\def\hend{$\glet\maybepush\relax \CDstrut\egroup}
\def\vbegin{\setbox\tempbox=\hbox\bgroup$}
\def\vend{$\egroup\ht\tempbox=\z@\dp\tempbox\CDvarrlen
  \box\tempbox}
\def\setCD{\let\harrow=\setCDarrow
  \let\CDCR=\setCR 
  \row=\@ne \col=\@ne \hbegin}
\let\endsetCD=\hend 

\def\findarrwd{\@rrwd=\z@ \count@=\col \advance\count@ by\sp@ncnt
  \loop\ifnum\count@>\col \advance\count@ by -1
      \advance\@rrwd by\csname col\the\count@\endcsname\repeat}
\def\setCDarrow#1#2#3{\kern\CDarrsurr\advance\col by \@ne
  \findarrwd \advance\@rrwd by -\r@wc@l  
  \@rrdp=\z@ 
  \maybepush
  \advance\col by -\@ne \advance\col by \sp@ncnt \span@ne
  \hbox to \@rrwd{\options
   \@rrwd=\scalefactor\@rrwd\hss
   \hplace{#1}\ulabel{#2}\dlabel{#3}\hss}%
   \kern\CDarrsurr}

\newdimen\labspacei 
\newdimen\labspaceii 

\newdimen\@xisheight
  \@xisheight=\the\fontdimen22\textfont2
\newdimen\labelskip
  \labelskip=\the\fontdimen10\textfont3 
\newdimen\topofshaft
\newdimen\botofshaft
\newdimen\botofulabel
\newdimen\topofdlabel
\def\getlabeldims{
  \topofshaft=0.5\sh@ftdiam
  \botofshaft=\topofshaft
  \advance\topofshaft by \@xisheight  
  \advance\botofshaft by -\@xisheight  
  \botofulabel=\topofshaft
  \advance\botofulabel by \labelskip
  \topofdlabel=\botofshaft
  \advance\topofdlabel by \labelskip}

\def\ulabel{\ifnum\row=\@ne\let\next\ulabeli
   \else\let\next\ulabellap\fi\next}
\def\ulabeli#1{\vbox{
  \clap{\kern-\@rrwd$\scriptstyle#1$}%
  \kern\botofulabel}\maybeoffset}
\def\ulabellap#1{\vbox to \z@{\vss
  \clap{\kern-\@rrwd$\scriptstyle#1$}%
  \kern\botofulabel}\maybeoffset}
\def\dlabel{\ifnum\row=\numrows\let\next\dlabeli
   \else\let\next\dlabellap\fi\next}
\def\dlabeli#1{\vtop{\kern\topofdlabel
  \clap{\kern-\@rrwd$\scriptstyle#1$}%
  }\maybeoffset}
\def\dlabellap#1{\vbox to \z@{\kern\topofdlabel
  \clap{\kern-\@rrwd$\scriptstyle#1$}%
  \vss}\maybeoffset}
\def\rlabel#1{\vbox to \z@{\vss
  \rlap{\kern\labelskip$\scriptstyle#1$}%
  \vss\kern-\@rrdp}\maybeoffset}
\def\llabel#1{\vbox to \z@{\vss
  \llap{$\scriptstyle#1$\kern\labelskip}%
  \vss\kern-\@rrdp}\maybeoffset}
\def\swlabel#1{\vtop{\kern0.5\@rrdp
  \llap{$\scriptstyle#1$\kern\labelskip\kern-0.5\@rrwd}
  }\maybeoffset}
\def\nwlabel#1{\vbox{
  \llap{$\scriptstyle#1$\kern\labelskip\kern-0.5\@rrwd}%
  \kern-0.5\@rrdp}\maybeoffset}
\def\selabel#1{\vtop{\kern0.5\@rrdp
  \rlap{\kern0.5\@rrwd\kern\labelskip$\scriptstyle#1$}%
  }\maybeoffset}
\def\nelabel#1{\vbox{
  \rlap{\kern0.5\@rrwd\kern\labelskip$\scriptstyle#1$}%
  \kern-0.5\@rrdp}\maybeoffset}
\def\cplace#1{\vbox to \z@{\vss
  \clap{$#1$\kern-\@rrwd}%
  \kern-\@rrdp\vss}\maybeoffset}
\def\hplace#1{\hbox to \@rrwd{#1}\maybeoffset}
\def\vplace#1{\clap{\vbox to \z@{#1\kern-\@rrdp}}\maybeoffset}

\newdimen\nudgeamount \nudgeamount=\z@
\newdimen\liftamount \liftamount=\z@
\let\maybeoffset\relax
\newbox\offsetbox \newdimen\lastheight
\def\dooffset{
  \setbox\offsetbox=\lastbox \lastheight=\ht\offsetbox 
  \setbox\offsetbox=\vbox{\kern-\liftamount\box\offsetbox}%
  \ht\offsetbox=\lastheight
  \kern\nudgeamount\box\offsetbox\kern-\nudgeamount
  \global\nudgeamount=\z@ \global\liftamount=\z@
  \glet\maybeoffset=\relax}
\def\nudge#1{\ifinmeasureCD\else
  \global\advance\nudgeamount#1\relax
  \global\let\maybeoffset\dooffset\fi}
\def\lift#1{\ifinmeasureCD\else
  \global\advance\liftamount#1\relax
  \global\let\maybeoffset\dooffset\fi}

\def\findarrdp{\@rrdp=\CDvarrlen
  \ifnum\sp@ncnt@>1
    \advance\@rrdp by \CDstrutlen
    \multiply\@rrdp by \sp@ncnt@
    \advance\@rrdp by -\CDstrutlen \fi
 }

\def\varrow#1#2#3{\ifnum\sp@ncnt>\@ne 
     \sp@ncnt@=\sp@ncnt\relax\fi
  \findarrdp \@rrwd=\z@ 
  \kern\c@l
   \hbox to \z@{\options
   \@rrdp=\scalefactor\@rrdp
    \hss\vplace{#1}\llabel{#2}\rlabel{#3}\hss}%
  \global\advance\col by \@ne \setsp@n\@ne\@ne
  }

\def\novarrow{\varrow\vfill{}{}}

\def\tweenarrows#1{\findarrwd \findarrdp \setsp@n\@ne\@ne
  \rlap{\options\cplace{#1}}}

\def\usarrow #1#2#3{\dr@tn=\@ne
  \findarrwd \findarrdp \setsp@n\@ne\@ne 
  \rlap{\options\cplace{#1}\nwlabel{#2}\selabel{#3}}%
  \dr@tn=\z@}
\def\dsarrow #1#2#3{\dr@tn=\tw@
  \findarrwd \findarrdp \setsp@n\@ne\@ne 
  \rlap{\options\cplace{#1}\swlabel{#2}\nelabel{#3}}%
  \dr@tn=\z@}
 \def\@rrow#1{\csname #1@rrow\endcsname}
 \def\R@rrow{\harrow \rtarrfill}
 \def\L@rrow{\harrow \ltarrfill}
 \def\V@rrow{\varrow \dnarrfill}
 \def\A@rrow{\varrow \uparrfill}
 \def\SE@rrow{\dsarrow \SEarrow}
 \def\NW@rrow{\dsarrow \NWarrow}
 \def\SW@rrow{\usarrow \SWarrow}
 \def\NE@rrow{\usarrow \NEarrow}
 \def\DS@rrow{\dsarrow \dnslope}
 \def\US@rrow{\usarrow \upslope}
 \def\upslope{\find@xyargs
       \@linelen=\unitlength\@sline}
 \def\dnslope{\find@xyargs\@yarg=-\@yarg\relax
       \@linelen=\unitlength\@sline}

\newtoks\optionlist 
\optionlist={}
\let\options\relax
\def\dooptions{\the\optionlist\global\optionlist={}%
  \glet\options=\relax}
\def\option#1{\ifinmeasureCD\else
  \glet\options=\dooptions
  \global\optionlist=\expandafter{\the\optionlist\relax#1}\fi}
\def\wider#1{\ifinmeasureCD\else
   \option{\advance\@rrwd by #1}\fi}
\def\deeper#1{\ifinmeasureCD\else
   \option{\advance\@rrdp by #1}\fi}

{\def\\{\global\let\sptoken= }\\ }

\def\CR{\futurelet\nexttok\testCR}
\def\testCR{\ifx\nexttok\sptoken
   \let\next\eatspaceCR\else\let\next\CDCR\fi\next}
\def\eatspaceCR#1 {\CR}
\def\measureCR{\ifx\nexttok\endmeasure\let\nextCR\relax
    \else\let\nextCR\drop\fi\nextCR}
\def\setCR{\ifodd\row
  \ifx\nexttok\endsetCD\else\hend\betweenCDrows\vbegin\fi
  \else\vend\betweenCDrows\hbegin\fi}

\countdef\dimenc@unt=11
\def\CD#1\endCD{
   \begingroup\let\\=\CR
  \m@th\offinterlineskip
   \measure#1\endmeasure\null\,\vcenter{\setCD#1\endsetCD}\,
   \endgroup
    }

\ifx\@clnwd\undefined \nol@g\else\catcode`\ =14\relax\fi
 \font\@linefnt=line10 
 \newcount\@tempcnta
 \newcount\@tempcntb
 \newdimen\@tempdima
 \newdimen\@tempdimb
 \newdimen\@wholewidth
 \newdimen\@halfwidth
   \@wholewidth\fontdimen8\@linefnt \@halfwidth .5\@wholewidth
 \newdimen\unitlength
 \newcount\@xarg
 \newcount\@yarg
 \newcount\@yyarg
 \newbox\@linechar
 \newdimen\@linelen
 \newdimen\@clnwd
 \newdimen\@clnht
 \newif\if@negarg
 
 \def\@whilenoop#1{}

 \def\@whiledim#1\do #2{\ifdim #1\relax#2\@iwhiledim{#1\relax#2}\fi}

 \def\@iwhiledim#1{\ifdim #1\let\@nextwhile=\@iwhiledim 
         \else\let\@nextwhile=\@whilenoop\fi\@nextwhile{#1}}

 \def\@sline{\ifnum\@xarg< 0 \@negargtrue \@xarg -\@xarg \@yyarg -\@yarg
   \else \@negargfalse \@yyarg \@yarg \fi
 \ifnum \@yyarg >0 \@tempcnta\@yyarg \else \@tempcnta -\@yyarg \fi
 \ifnum\@tempcnta>6 \@badlinearg\@tempcnta0 \fi
 \ifnum\@xarg>6 \@badlinearg\@xarg 1 \fi
 \setbox\@linechar\hbox{\@linefnt\@getlinechar(\@xarg,\@yyarg)}%
 \ifnum \@yarg >0 \let\@upordown\raise \@clnht\z@
    \else\let\@upordown\lower \@clnht \ht\@linechar\fi
 \@clnwd=\wd\@linechar
 \if@negarg \hskip -\wd\@linechar \def\@tempa{\hskip -2\wd\@linechar}\else
      \let\@tempa\relax \fi
 \@whiledim \@clnwd <\@linelen \do
   {\@upordown\@clnht\copy\@linechar
    \@tempa
    \advance\@clnht \ht\@linechar
    \advance\@clnwd \wd\@linechar}%
 \advance\@clnht -\ht\@linechar
 \advance\@clnwd -\wd\@linechar
 \@tempdima\@linelen\advance\@tempdima -\@clnwd
 \@tempdimb\@tempdima\advance\@tempdimb -\wd\@linechar
 \if@negarg \hskip -\@tempdimb \else \hskip \@tempdimb \fi
 \multiply\@tempdima \@m
 \@tempcnta \@tempdima \@tempdima \wd\@linechar \divide\@tempcnta \@tempdima
 \@tempdima \ht\@linechar \multiply\@tempdima \@tempcnta
 \divide\@tempdima \@m
 \advance\@clnht \@tempdima
 \ifdim \@linelen <\wd\@linechar
    \hskip \wd\@linechar
   \else\@upordown\@clnht\copy\@linechar\fi}
 
 \def\@getlinechar(#1,#2){\@tempcnta#1\relax\multiply\@tempcnta 8
 \advance\@tempcnta -9 \ifnum #2>0 \advance\@tempcnta #2\relax\else
 \advance\@tempcnta -#2\relax\advance\@tempcnta 64 \fi
 \char\@tempcnta}
 
 \def\vector(#1,#2)#3{\@xarg #1\relax \@yarg #2\relax
 \@tempcnta \ifnum\@xarg<0 -\@xarg\else\@xarg\fi
 \ifnum\@tempcnta<5\relax
 \@linelen=#3\unitlength
 \ifnum\@xarg =0 \@vvector 
   \else \ifnum\@yarg =0 \@hvector \else \@svector\fi
 \fi
 \else\@badlinearg\fi}
 
 \def\@svector{\@sline
 \@tempcnta\@yarg \ifnum\@tempcnta <0 \@tempcnta=-\@tempcnta\fi
 \ifnum\@tempcnta <5
   \hskip -\wd\@linechar
   \@upordown\@clnht \hbox{\@linefnt  \if@negarg 
   \@getlarrow(\@xarg,\@yyarg) \else \@getrarrow(\@xarg,\@yyarg) \fi}%
 \else\@badlinearg\fi}
 
 \def\@getlarrow(#1,#2){\ifnum #2 =\z@ \@tempcnta='33\else
 \@tempcnta=#1\relax\multiply\@tempcnta \sixt@@n \advance\@tempcnta
 -9 \@tempcntb=#2\relax\multiply\@tempcntb \tw@
 \ifnum \@tempcntb >0 \advance\@tempcnta \@tempcntb\relax
 \else\advance\@tempcnta -\@tempcntb\advance\@tempcnta 64
 \fi\fi\char\@tempcnta}
 
 \def\@getrarrow(#1,#2){\@tempcntb=#2\relax
 \ifnum\@tempcntb < 0 \@tempcntb=-\@tempcntb\relax\fi
 \ifcase \@tempcntb\relax \@tempcnta='55 \or 
 \ifnum #1<3 \@tempcnta=#1\relax\multiply\@tempcnta
 24 \advance\@tempcnta -6 \else \ifnum #1=3 \@tempcnta=49
 \else\@tempcnta=58 \fi\fi\or 
 \ifnum #1<3 \@tempcnta=#1\relax\multiply\@tempcnta
 24 \advance\@tempcnta -3 \else \@tempcnta=51\fi\or 
 \@tempcnta=#1\relax\multiply\@tempcnta
 \sixt@@n \advance\@tempcnta -\tw@ \else
 \@tempcnta=#1\relax\multiply\@tempcnta
 \sixt@@n \advance\@tempcnta 7 \fi\ifnum #2<0 \advance\@tempcnta 64 \fi
 \char\@tempcnta}
\catcode`\ =10

\dol@g 
\catcode`\@=\active
\LaTeXfeatures

\newtheorem{definition}{Definition}[section]
\newtheorem{theorem}{Theorem}
\newtheorem{lemma}[definition]{Lemma}
\newtheorem{proposition}{Proposition}
\newtheorem{corollary}[definition]{Corollary}
\newtheorem{remark}[definition]{Remark}
\newtheorem{problem}[definition]{Problem}

\newenvironment{proof}{\hskip-\parindent{\sc Proof}.\ \ }{
                       \hfill$\Box$\vskip\partopsep \vskip\topsep}

\def\defemb#1#2{\expandafter\def\csname #1\endcsname
                              {\relax\ifmmode #2\else\hbox{$#2$}\fi}}

\defemb{cB}{{\cal B}}
\defemb{cC}{{\cal C}}
\defemb{cF}{{\cal F}}
\defemb{cG}{{\cal G}}
\defemb{cI}{{\cal I}}
\defemb{cL}{{\cal L}}
\defemb{cU}{{\cal U}}
\defemb{cK}{{\cal K}}

\newcommand{\ct}{\centerline}

 \newcommand{\Aut}{\mathrm{Aut}}

\newcommand{\Int}{\mathrm{Int}}
\newcommand{\End}{\mathrm{End}}
\newcommand{\Ker}{\mathrm{Ker}}

\newcommand{\GL}{\mathrm{GL}}
\newcommand{\Hom}{\mathrm{Hom}}
\newcommand{\Var}{\mathrm{Var}}
\newcommand{\Pol}{\mathrm{Pol}}

\newcommand\K{\mathrm{K}}

\begin{document}
   \baselineskip 18pt



\ct{\Large Automorphisms of the category of free Lie algebras}

\bigskip

\ct{\Large G.~Mashevitzky ${}^{\natural}$, B.~Plotkin ${}^{\sharp}$,
  E.Plotkin ${}^{\flat}$\footnote{This work is partially supported by the
  Emmy Noether Institute for
Mathematics and the Minerva Foundation of Germany, by the
Excellency Center of the Israel science Foundation "Group
Theoretic Methods in the Study of Algebraic Varieties" and by
INTAS 00-566}}

\bigskip

 \ct {${}^{\natural}$ \Large{\it {Ben Gurion University of the Negev}}}

 \bigskip

\ct{ ${}^{\sharp}$ \Large{\it {Hebrew University}}}

\bigskip

\ct {${}^{\flat}$\Large {\it { Bar Ilan University}}}




\begin{abstract}
We prove that every automorphism of the category of free Lie
algebras is a semi-inner automorphism. This solves the problem 3.9
from \cite{ MPP} for Lie algebras.
\end{abstract}

\section*{Introduction}
We start from an arbitrary variety of algebras $\Theta$. Let us
denote the category of free in $\Theta$ algebras $F=F(X)$, where $X$
is finite, by $\Theta^0$. In order to avoid the set
theoretic problems we view all $X$ as subsets of a universal
infinite set $X^0$.

Our main goal is to study automorphisms $\varphi: \Theta^0\to
\Theta^0$ and the corresponding group $\Aut \Theta^0$ for various
$\Theta$.

In this paper we consider the case when $\Theta$ is the variety of
all Lie algebras over an infinite field $P$. Our aim is to  prove
the following principal theorem:

\begin{theorem}
Every automorphism of the category of free Lie algebras
 is a semi-inner automorphism.
\end{theorem}

  This Theorem solves the problem 3.9 from \cite{MPP} for the case of
 Lie algebras.

Our primary interest to automorphisms of  categories raised from
the universal algebraic geometry (see \cite{Pl5}, \cite{Pl4},
\cite{Pl2}, \cite{Pl3}, \cite{NP}, \cite{Be1},
\cite{BMR1},\cite{BMR2},\cite{MR1},\cite{MR2},\cite{Se}, etc). The
motivations we keep in mind are inspired by the following
observations.

Some basic notions of classical algebraic geometry can be defined
for arbitrary varieties of algebras $\Theta$. For every algebra
$H\in \Theta$ one can consider geometry in $\Theta$ over $H$. This
geometry gives rise to the category $\K_\Theta(H)$ of algebraic
sets
 in affine spaces over $H$ \cite{Pl5}. The key question
in this setting is  when the geometries in $\Theta$ defined by
different algebras $H_1$ and $H_2$ coincide. The coincidence of
geometries means for us that the corresponding categories of
algebraic sets $\K_\Theta(H_1)$ and $\K_\Theta(H_2)$ are either
isomorphic or equivalent.

It is known that the conditions on $H_1$ and $H_2$ providing
isomorphism or equivalence of the categories $\K_\Theta(H_1)$ and
$\K_\Theta(H_2)$ depend essentially on the description of the
automorphisms of the category $\Theta^0$ (see \cite{Pl5},
\cite{MPP}). This explains the interest to automorphisms of
categories of free algebras of varieties.

Let $F=F(X)\in\Theta$ be a free algebra, i.e., an object of the
category $\Theta^0$. The group $\Aut(\Theta^0)$ is tied naturally
with
the following sequence of groups:
$$
\Aut(F), \Aut(\Aut(F)), \Aut(End(F)) 
$$
The groups $\Aut(F)$ are known for the variety of all groups
(Nielsen's theorem, \cite{LS}),
for the variety of Lie algebras (P.Cohn's theorem, \cite{C}),
for the free associative algebras over a field when
the number of generators of $F$ is $\leq 2$  (\cite {Co},
\cite{LML},\cite{Cz}, \cite{Na}) and for some other varieties. For
free associative algebras with bigger number of generators the
question is still open (see Cohn's conjecture \cite{Co}). The
groups $\Aut(\Aut(F))$, $\Aut(End(F))$ are known for the variety
of all groups \cite{DF}, \cite{For}, \cite{T1},
and due to E.Formanek every automorphism of $\End(F)$  is inner.
The groups $\Aut(\Aut(F))$, $\Aut(End(F))$ are also known for some other
varieties of groups and semigroups \cite{MS}, \cite{ST}, \cite{SH},
\cite{DF1}, \cite{T2}.

Suppose that a free algebra $F=F(X)$ generates the whole variety
$\Theta$. In this case there exists a natural way from the group
$\Aut(End(F))$ to the group $\Aut(\Theta^0)$. Thus, there is a
good chance  to reduce the question on automorphisms of the
category $\Theta^0$ to the description of $\Aut(End(F))$.

$\Aut (F)$ is the group of  invertible elements of the
semigroup $\End (F)$. Every automorphism $\varphi$ of the
semigroup $\End (F)$ induces an automorphism of the group $\Aut
(F)$. This gives a homomorphism $\tau:\Aut(\End(F))\to
\Aut(\Aut(F))$.  The kernel of this homomorphism consists of
automorphisms acting trivially in $\Aut (F)$. These automorphisms
are called {\it stable}. We will prove that
\begin{enumerate}\item
The homomorphism $\tau$ is not a surjection.
\item
If $X$ consists of more than 2 elements then $\tau$ is an
injection.
\item
If $X$ consists of 2 elements then $\Ker \tau$ consists of
scalar automorphisms (see Sections \ref{2}, \ref{3}).
\end{enumerate}

The paper is organized as follows. In  Section 1 we give the
definitions of inner and semi-inner automorphisms of a
category. In Section 2 the notations are introduced.
Section 3 is dedicated to linearly stable
automorphisms and we prove that every linearly stable automorphism
is inner. In Section 4 we define the notion of a quasi-stable automorphism
and prove that every quasi-stable automorphism is inner.
 In Section 5 we prove that every
automorphism of the semigroup of endomorphisms of the free two
generator Lie algebra is semi-inner. In Section 6 we prove the
general reduction theorem for a large class of varieties
and reduce the
problem about $\Aut (\Theta^0)$ to the description of
$\Aut(End(F(x,y))$. Section 7 is dedicated to the proof of the main
theorem. Finally, in Appendix we prove some auxiliary statements
used in the text.

\section{Inner and semi-inner automorphisms of a category}

Recall the notions of category isomorphism and equivalence
\cite{ML}. A functor $\varphi:{\cal C}_1 \to {\cal C}_2$ is called
an {\it isomorphism of categories} if there exists a functor
$\psi: {\cal C}_2 \to
{\cal C}_1$ such that $\psi\varphi=1_{{\cal C}_1}$ and
$\varphi\psi=1_{{\cal C}_2}$, where $1_{{\cal C}_1}$ and $1_{{\cal
C}_2}$ are identity functors.

Let $\varphi_1$, $\varphi_2$ be two functors ${\cal C}_1
\to {\cal C}_2$. An {\it isomorphism of functors}
 $s: \varphi_1\to\varphi_2$ is defined by the following conditions:\\
 1. To every object $A$ of the category $C_1$ an isomorphism $s_A:\varphi_1(A)\to\varphi_2(A)$ in
${\cal C}_2$ is assigned.\\
2. If $\nu:A\to B$ is a morphism in ${{\cal C}_1}$ then there is a
 commutative diagram in ${\cal C}_2$:

$$\CD \varphi_1(A) @>s_A>>\varphi_2(A) \\ @V\varphi_1(\nu)VV @VV\varphi_2(\nu)
 V\\\varphi_1(B)@>s_B>>\varphi_2(B) \endCD$$

The isomorphism of functors $\varphi_1$ and $\varphi_2$ is denoted by
 $\varphi_1\simeq\varphi_2$.

The notion of category equivalence generalizes the notion of
category isomorphism. A pair of functors $\varphi: {\cal C}_1\to
{\cal C}_2$ and $\psi: {\cal C}_2 \to{\cal C}_1$ define a {\it
category equivalence} if $\psi\varphi\simeq 1_{{\cal C}_1}$ and
$\varphi\psi\simeq 1_{{\cal C}_2}$. If ${\cal C}_1= {\cal
C}_2={\cal C}$ then we get the notions of {\it automorphism} and
{\it autoequivalence} of the category ${\cal C}$.
\begin{definition}
 An automorphism $\varphi$ of the category $\cal C$
is called inner if there exists an isomorphism of functors
$s:1_{\cal C}\to \varphi$.
\end{definition}

This means that for every object $A$ of the category $\cal C$
there exists an isomorphism $s_A: A\to \varphi(A)$ such that
$$
\varphi(\nu)=s_B\nu s_A^{-1}: \varphi(A)\to \varphi(B),
$$
for any morphism $\nu: A\to B$ in $\cal C$.

For every small category $\cal C$ denote the group of all automorphisms
of $\cal C$ by $\Aut(\cal C)$ and denote its
normal subgroup of all inner automorphisms by $\Int (\cal C)$.

From now on and till  Section 5 $\Theta$ will denote the variety
of all Lie algebras over the field $P$. Correspondingly,
$\Theta^0$ is the category of free Lie algebras over $P$.

Define the notion of a semi-inner automorphism of the category
$\Theta^0$. Consider, first, semimorphisms in the variety $\Theta$.
A {\it semimorphism} in $\Theta$ is a pair
 $(\sigma,\nu): A\to B$, where $A$ and $B$ are algebras in $\Theta$,
$\nu: A\to B$ a homomorphism of Lie rings, $\sigma$ an
automorphism of the field $P,$ subject to condition $\nu(\lambda
a)=\sigma(\lambda)\nu(a)$, where $\lambda\in P$, $a\in A$. If
$\sigma =1$ then $\nu$ is a homomorphism of Lie algebras and we
write it as $(1,\nu)$. Semimorphisms are multiplied componentwise.

Thus, if $\mu: A\to B$
is a homomorphism, and $(\sigma, \nu_1): A\to A_1$, $(\sigma,
\nu_2): B\to B_1$ are semi-isomorphisms, then
$$
(\sigma, \nu_2)(1,\mu)(\sigma, \nu_1)^{-1}= ((\sigma,
\nu_2)(1,\mu)(\sigma^{-1}, \nu_1^{-1})= (1, \nu_2\mu\nu_1^{-1}).
$$
This means that $\nu_2\mu\nu_1^{-1}: A_1\to B_1$ is a homomorphism
of Lie algebras.

\begin{definition} An automorphism $\varphi$ of the category $\Theta^0$ is
called semi-inner, if for some $\sigma\in \Aut (P)$ there is a
semi-isomorphism of functors $(\sigma, s):1_{\Theta^0}\to \varphi$.
\end{definition}

This definition means that for every object $F\in \Theta^0$ there
exists  a semi-isomorphism $(\sigma, s_F): F\to \varphi(F)$  such
that
$$
\varphi(\nu)=s_{F_2}\nu s_{F_1}^{-1}:\
\varphi(F_1)\to\varphi(F_2),
$$
for any morphism $\nu : F_1\to F_2$ in $\Theta^0$.

Let now $\sigma$ be an arbitrary automorphism of the field $P$. We
will construct a semi-inner automorphism $\hat{\sigma}$ of the
category $\Theta^0$. Consider an arbitrary free finitely generated
Lie algebra $F=F(X)$ and fix a Hall basis in $F$
\cite{Ba},\cite{BK}. It can be shown \cite{BK} that if $u,v$ are
two elements of a  Hall
basis then $[u,v]$ is presented via elements of this basis with the
coefficients (structure constants) belonging to the minimal
subfield of the field $P$. Define a map $\sigma_F:F \to F$. Every
element $w$ of $F$ has the form:
$$
w=\lambda_1 u_1+\cdots+\lambda_n u_n,
$$
where $\lambda_i\in P$ and $u_i$ belong to the Hall basis of $F$.
Define
$$
\sigma_F(w)=\sigma(\lambda_1) u_1+\cdots+\sigma(\lambda_n) u_n.
$$

We show that $(\sigma,\sigma_F)$ is a semi-automorphism of the
algebra $F$. It is clear that $\sigma_F$ preserves the addition, and
that $\sigma_F(\lambda w)=\sigma(\lambda)\sigma_F(w)$. It remains
to check that $\sigma_F$ preserves the multiplication. Let
$w_1=\sum_i\alpha_i u_i$, and $w_2=\sum_j\alpha_j u_j$. Take
$[u_i,u_j]=\sum_k n_k^{i,j}u_k^{i,j}$, where $n_k^{i,j}$ belong to
a minimal subfield . Then
$$
[w_1,w_2]=\sum_{i,j}\alpha_i\beta_j[u_i,u_j]=\sum_{i,j,k}\alpha_i\beta_j
 n_k^{i,j}u_k^{i,j}.
$$
Apply $\sigma_F$, then
$$
\sigma_F[w_1,w_2]=\sum_{i,j,k}\sigma(\alpha_i)\sigma(\beta_j)
n_k^{i,j}u_k^{i,j},
$$
since the automorphism $\sigma $ does not change the elements of a
prime subfield. On the other hand
$$
[\sigma_F(w_1),\sigma_F(w_2)]=[\sum_i\sigma(\alpha_i)u_i,
\sum_j\sigma(\beta_j)u_j]=
\sum_{i,j}\sigma(\alpha_i)\sigma(\beta_j)[u_i,u_j]=
\sum_{i,j,k}\sigma(\alpha_i)\sigma(\beta_j) n_k^{i,j}u_k^{i,j}.
$$
We verified that the pair $(\sigma, \sigma_F)$ defines a
semi-automorphism of the algebra $F$. Note that $\sigma_F$ does
not change variables from $X$ and does not change all commutators
constructed from these variables.

Now we are able to define the automorphism $\hat\sigma$ of the
category $\Theta^0$. This automorphism does not change objects and
$\hat\sigma(\nu)=\sigma_{F_2}\nu\sigma_{F_1}^{-1}: F_1 \to F_2$
for every morphism $\nu: F_1\to F_2$.
It is easy to check that if $\varphi$ is an arbitrary
semi-inner automorphism of $\Theta^0$ for the given $\sigma\in\Aut
(P)$ then there is a factorisation $\varphi=\varphi_0\hat\sigma$,
where $\varphi_0$ is an inner automorphism.

The same scheme which was applied for the definition of semi-inner
automorphisms of the category $\Theta^0$ works for the definition
of semi-inner automorphisms of the semigroup $\End(F)$, where $F$
is a free finitely generated Lie algebra. An automorphism
$\varphi$ of $\End (F)$ is called a {\it semi-inner} automorphism
if  there exists a semi-automorphism $(\sigma, s):F \to F$ such
that $\varphi(\nu)=s\nu s^{-1}$ for every $\nu \in \End (F)$. The
factorisation $\varphi=\varphi_0\hat\sigma$, where $\varphi_0$ is
inner holds also in this case.

\section{Notations and preliminaries}\label{2}

Let $X$ be a finite set. Denote by $F(X)$ the free Lie algebra
over an infinite field $P$
generated by the set $X$ of free generators. The Lie operation is denoted
by $[,]$. Denote the group of all non zero elements of $P$ by $P^*$.
We denote the semigroup of all
endomorphisms of $F(X)$ by $\End(F(X))$. Any endomorphism of
$F(X)$ is uniquely determined by a mapping $X\rightarrow F(X)$.
Therefore, we define an endomorphism $\varphi$ of $F(X)$ by
defining $\varphi(x)$ for all $x\in X$.
 Denote the group of all automorphisms of $F(X)$ by $\Aut(F(X))$.

We fix a basis, say the Hall basis in $F(X)$ and consider the
presentations of elements of $F(X)$ in this basis.
Denote the length of a monomial $u\in F(X)$ by $|u|$, we call it
also the degree of $u$.
Denote the set of all elements of $X$ included in $u$ by
$\chi(u)$. The set $\chi(u)$ is a support of element $u$, which
is uniquely defined by the presentation of $u$ in the fixed basis.
Denote the number of occurrences of a letter $x$ in $u$ by
$l_x(u)$. Let $p\in F(X)$. Denote the degree of the polynomial $p$
by $deg(p)$. Denote the cardinality of the set $X$ by $|X|$.

Let us denote the semigroup of all endomorphisms $\varphi\in
\End(F(X))$ which assign a linear polynomial from $F(X)$ to any $x\in X$
by ${\End}_l(F(X))$. Denote the group of the linear automorphisms of
$F(X)$ by ${\Aut}_l(F(X))$.

Let $X=\{ x_1\ldots x_n\}$. If $\varphi$ is linear then
$\varphi(x_i)= a_{i1}x_1+\cdots +a_{in}x_n$, where $a_{ij}\in P$.
This means that a linear automorphism is defined by its matrix of
coefficients. The multiplication of linear automorphisms corresponds
to the multiplication of
 their matrices. Thus, the semigroup ${\End}_l(F(X))$ is isomorphic to the
matrix semigroup $M_n(P)$. The scalar matrix corresponds to the
linear automorphism, defined by $f_a(x_i)=ax_i$. Therefore,
automorphism $f_a$ commutes with all linear endomorphisms.
However, it does not commute with an arbitrary endomorphism.

Denote the endomorphism of $F(X)$ which assigns the same $p\in
F(X)$ to any $x\in X$ by $c_p$. All endomorphisms of the type
$c_p$ form a subsemigroup $C_F$ of $\End(F(X))$. Let us denote the
subsemigroup of $C_F$ consisting of all endomorphisms of the type
$c_u$ where $u\in X$ by $C_X$. Let $C_l=C_F\cap {\End}_l(F(X))$.
$C_l$ consists of all endomorphisms of the type $c_p\in C_F$ where
$p$ is a linear polynomial.

\begin{definition}
An automorphism $\xi \in \Aut(\End(F(X)))$ which acts identically
on ${\Aut}_l(F(X))$, is called a linearly stable automorphism of
$\End(F(X))$.
\end{definition}


\section{Linearly stable automorphisms of $\End(F(X))$ are inner}\label{3}

In this section we prove that the semigroup $C_F$ of all constant
endomorphisms is invariant in respect to the action of any
linearly stable automorphism of $\End(F(X))$ and that any linearly
stable automorphism acts identically on the semigroup
${\End}_l(F(X))$ of all linear endomorphisms. Then we prove that
any linearly stable automorphism of $\End(F(X))$ is an inner
automorphism.

\begin{lemma}\label{C_F}
Let $\xi$ be a linearly stable automorphism of $\End(F(X))$. Then
$\xi(C_F)=C_F$.
\end{lemma}

\begin{proof}
Take $u\in F(X)$, $c_u\in C_F$.
Consider $g\in\Aut_l(F(X))$ such that $g(x)=y, g(y)=x,$ and
$g(z)=z$ for any  $z\in X$ distinct from $x$ and $y$. Then
$c_ug=c_u$ and $\xi(g)=g$. Hence, $\xi(c_u)(x)=\xi(c_ug)(x)=
\xi(c_u)\xi(g)(x)=\xi(c_u)(y)=v$. Thus, $\xi(c_u)=c_v \in C_F$.
Therefore, $\xi(C_F)\subset C_F$. Similarly, $\xi^{-1}(C_F)\subset
C_F$. Hence, $\xi(C_F)=C_F$.
\end{proof}

\begin{lemma}\label{homogenity}
Let $\xi\in \Aut(\End(F(X)))$ be a linearly stable automorphism.
Let $\xi(c_p)=c_q$. Then $\chi(p)=\chi(q)$ (polynomials $p$ and
$q$ are constructed from the same variables).
\end{lemma}

\begin{proof}
Suppose that $x\in\chi(p)\setminus\chi(q)$. Let $a\in P$ be an
element of infinite order (torsion free). Consider $g_a\in
{\Aut}_l(F(X))$, which assigns $ax$ to $x$ and assigns $z$ to $z$
for any $z\in X$ different from $x$. Then $g_a(p)=w\neq p$ (indeed
all elements of the basis are linearly independent and $a$ is
torsion free) and $g_a(q)=q$. Therefore, $g_ac_p=c_w\neq c_p$. On
the other hand,
$\xi(g_ac_p)=\xi(g_a)\xi(c_p)=g_ac_q=c_q=\xi(c_p)$. It contradicts
to injectivity of $\xi$. The similar reasoning can be applied to
$x\in \chi(q)\setminus\chi(p)$ as well. Thus, $\chi(q)=\chi(p)$.
\end{proof}

\begin{lemma}\label{C_X}
Any linearly stable automorphism $\xi$ of $\End(F(X))$
acts identically on $C_X$.
\end{lemma}

\begin{proof}
It is obvious that any $c_x\in C_X$ is a right identity of $C_F$.
Therefore, $\xi(c_x)=c_p\in C_F$ (Lemma \ref{C_F}) and $\xi(c_x)$
is a right identity of $C_F$. Hence, $c_xc_p=c_x$. It follows from
Lemma \ref{homogenity} and identity $[xx]=0$ that $p=ax$, $a\in
P$. Therefore, $c_xc_p=c_p$. Thus, $\xi(c_x)=c_p=c_x$.
\end{proof}

\begin{lemma}\label{C-def}
For any linearly stable automorphism $\xi$ of $\End(F(X))$ and for any
$\varphi \in \End (F(X))$ $\xi(c_{\varphi(x)})=c_{\xi(\varphi)(x)}$.
\end{lemma}

\begin{proof}
We have $\varphi c_x=c_{\varphi(x)}$. Then
$\xi(c_{\varphi(x)})=\xi(\varphi c_x)=\xi(\varphi)
c_x=c_{\xi(\varphi)(x)}$.
\end{proof}

\begin{lemma}\label{E}
 Any linearly stable automorphism $\xi$ of
$\End(F(X))$ acts identically on ${\End}_l(F(X))$.
\end{lemma}

\begin{proof}
Let us prove first that $\xi$ acts identically on $C_l$. Let
$c_p\in C_l$, where $p$ is a linear polynomial. Let $g\in
{\Aut}_l(F(X))$ be defined by $g(x)=p$, $x\in \chi(p)$, and $g(y)=y$
for any other $y\in X$.  It follows from Lemma \ref{C_X} that
$\xi(c_x)=c_x$. Therefore, $gc_x=c_p$ and
$\xi(c_p)=\xi(g)\xi(c_x)=gc_x=c_p$.

Hence,
$c_{\varphi(x)}=\xi(c_{\varphi(x)})=\xi(\varphi)c_x=c_{\xi(\varphi)(x)},$
for any $\varphi\in End_l(F(X))$ . Thus,
$\xi(\varphi)(x)=\varphi(x)$.
\end{proof}

\begin{lemma}\label{lincomb}
Let $\xi\in \Aut(End(F(X)))$ be a linearly stable automorphism.
Let $\xi(c_{p_i})=c_{q_i}$ and $a_i\in P$ for $i=1,2,...,k$. Then
$\xi(c_{a_1p_1+...+a_kp_k})=c_{a_1q_1+...+a_kq_k}$.
\end{lemma}

\begin{proof}
It follows from Lemma \ref{E} that $\xi(c_{ax})=c_{ax}$.
Therefore,
$\xi(c_{ap_i})=\xi(c_{p_i}c_{ax})=c_{q_i}c_{ax}=c_{aq_i}$.

Define $\varphi\in \End(F(X))$ as follows $\varphi(x)=p_1,
\varphi(y)=p_2$.
Then $c_{p_1+p_2}=\varphi c_{x+y}$. It follows from Lemma
\ref{C-def} that $\xi(\varphi)(x)=q_1$ and $\xi(\varphi)(y)=q_2$.
It follows from Lemma \ref{E} that $\xi(c_{x+y})=c_{x+y}$.
Therefore, $\xi(c_{p_1+p_2})=\xi(\varphi) c_{x+y}=c_{q_1+q_2}$.

Let us prove that
$\xi(c_{a_1p_1+...+a_kp_k})=c_{a_1q_1+...+a_kq_k}$ by  induction
on $k$. We have proved above the basis of the induction:
$\xi(c_{a_1p_1})=c_{a_1q_1}$. Assume that the statement is true
for $k<t$.
It follows from the assumption of the induction that
$\xi(c_{a_1p_1+...+a_{t-1}p_{t-1}})=c_{a_1q_1+...+a_{t-1}q_{t-1}}$
and $\xi(c_{a_tp_t})=c_{a_tq_t}$. $\xi(c_{p+p'})=c_{q+q'}$.
Therefore, $\xi(c_{a_1p_1+...+a_tp_t})=c_{a_1q_1+...+a_tq_t}$.
\end{proof}

\begin{lemma}\label{axy}
Let $\xi\in \Aut(\End(F(X)))$ be a linearly stable automorphism. Then
$\xi(c_{[x_1,x_2]})=c_{a[x_1,x_2]}$, where $a\in P$ and $a\neq 0$.
\end{lemma}

\begin{proof}
Suppose that $\xi(c_{[x_1,x_2]})=c_f$.

Denote the endomorphism generated by the mapping $x_i\rightarrow
a_ix_i$ by $\tau_{a_1...a_n}$. It follows from Lemma \ref{E} that
$\xi(\tau_{a_1...a_n})=\tau_{a_1...a_n}$. Therefore,
$\xi(\tau_{a_1...a_n}c_{[x_1,x_2]})=\tau_{a_1...a_n}c_f$. Hence,
$\xi(c_{a_1a_2[x_1,x_2]})=c_{\tau_{a_1...a_n}(f)}$. It follows from
Lemma \ref{lincomb} that $\xi(c_{a_1a_2[x_1,x_2]})=c_{a_1a_2f}$.
Thus, $a_1a_2f=\tau_{a_1...a_n}(f)$.

Suppose that $f=\alpha_1f_1+...+\alpha_tf_t$ is the decomposition
of $f$ with respect to the basis of Hall. Observe that the
decomposition of $\tau_{a_1...a_n}(f)$ contains the same elements
of the basis as the decomposition of $f$ but with different
coefficients. Elements of the basis are linearly independent over
$P$. Therefore, we obtain a system of equations of the form
$(a_1^{k^i_1}...a_n^{k^i_n}-a_1a_2)\alpha_i=0$ from the equation
$a_1a_2f=\tau_{a_1...a_n}(f)$. For the element $[x_1,x_2]$ of the
basis we get the equation $(a_1a_2-a_1a_2)\alpha=0$. In all other
cases it is easy to find $a_1,...,a_n$ such that
$(a_1^{k^i_1}...a_n^{k^i_n}-a_1a_2)\neq 0$ (for example, if
$char(P)\neq 2$, take $a_1=a_2=...=a_n=2$ for the monomials with
the number of multipliers different from $2$ and take $a_1=a_2=2,
a_3=...=a_n=1$ for monomials $[x_i,x_j]$, where
$\{i,j\}\neq\{1,2\}$). Hence, all coefficients $\alpha_i$ are
equal to zero except the coefficient of the monomial $[x_1,x_2]$.
\end{proof}


Let $a\in P^*$. Consider a scalar automorphism $f_a$ of
$F(X)$ defined by the rule $f_a(x)=ax$, for every $x\in X$. It
defines an inner automorphism ${\hat f}_a$ of the semigroup
$\End(F(X))$ by the rule ${\hat f}_a(\varphi)=f_a\varphi
f_a^{-1}$. ${\hat f}_a$ acts trivially on ${\End}_l(F(X))$. Thus,
${\hat f}_a$ is linearly stable.

\begin{proposition}\label{auto-a}
Let $\xi\in \Aut(\End(F(X)))$ be a linearly stable automorphism.
Then there exists $a\in K$ such that
 $\xi ={\hat f}_a$.
\end{proposition}

\begin{proof}
Let $\xi(c_{[x_1,x_2]})=c_{a[x_1,x_2]}$. Take the scalar
automorphism $f_a$ corresponding to the element $a$. Consider the
bijection $F(X)\rightarrow F(X)$ which multiplies a monomial of
the length $n$ on $a^{n-1}$. Suppose $p$ is a polynomial presented
as the sum of its homogeneous components: $p=p_1+...+p_s$, where
$deg(p_i)=i$ or $p_i=0$. Denote $p_1+ap_2+...+a^{s-1}p_s$ by
$\bar{p}$. Let us prove that $\xi(c_p)=c_{\bar {p}}$ for any
$c_p\in C_F$ by  induction on the number $r$ of monomials of the
polynomial $p$.

Let us prove the base of the induction for $r=1$ by induction on
the degree of the monomial $p=u$. The base of this induction
follows from the Lemma \ref{E}. Suppose that
$\xi(c_u)=c_{a^{l-1}u}$ for any $u$ such that $|u|=l<k$. Suppose
now that $|u|=k$ and $u=[u_1,u_2]$, where $|u_1|=k_1$ and
$|u_2|=k_2$. Let $\varphi(x_1)=u_1, \varphi(x_2)=u_2$ and
$\varphi(x)=x$ for any other $x\in X$. Then $\varphi
c_{[x_1,x_2]}=c_u$. It follows from the assumption of the
induction that $\xi(c_{u_1})=c_{a^{k_1-1}u_1}$,
 $\xi(c_{u_2})=c_{a^{k_2-1}u_2}$ and
$\xi(c_x)=c_x$ for any $x\in X$. Therefore, it follows from Lemma
\ref{C-def} that $\xi(\varphi)(x)=\overline{\varphi(x)}$ for any
$x\in X$. Hence, $\xi(c_u)=\xi(\varphi c_{[x_1,x_2]})=\xi(\varphi)
c_{a[x_1,x_2]}=
c_{\xi(\varphi)(a[x_1,x_2])}=c_{a[a^{k_1-1}u_1,a^{k_2-1}u_2]}=c_{a^{k-1}u}
=c_{\bar{u}}$. Thus, we proved the basis of the first induction
for $r=1$.

Suppose that $\xi(c_p)=c_{\bar {p}}$ for any $ p\in F(X)$ which
contains less than $k$ monomials. Suppose now that $p$
contains $k$ monomials and $p=q+g$ where each of the polynomials
$q$ and $g$ contains less than $k$ monomials. It follows from the
assumption of the induction that $\xi(c_{q})=c_{\bar {q}}$,
$\xi(c_{g})=c_{\bar {g}}$ and $\xi(c_x)=c_x$ for any $x\in X$.
Hence, it follows from Lemma \ref{lincomb} that
$\xi(c_p)=\xi(c_{q+g})=c_{\bar {q}+\bar {g}}=c_{\bar {p}}$. Thus,
$\xi(c_p)=c_{\bar {p}}$ for any $p\in F(X)$.

In particular,
$c_{\xi(\varphi)(x)}=\xi(\varphi c_x)=\xi(c_{\varphi(x)})=
c_{\overline{\varphi(x)}}$ for every $\varphi \in \End(F(X))$. Thus,
$\xi(\varphi)(x)=\overline{\varphi(x)}$. On the other hand it is easy
to see that ${\hat f}_a(\varphi)(x)=\overline{\varphi(x)}$. Therefore,
$\xi(\varphi)(x)=
 {\hat f}_a(\varphi)(x)$. This equality holds for every $x\in X$ and
every $\varphi\in \End(F(X))$. Thus, $\xi={\hat f}_a$.
\end{proof}

\begin{corollary}\label{lss}
Any linearly stable automorphism of $\End(F(X))$ is inner.
\end{corollary}

\begin{remark}\label{R}
1. Let $|X|>2$. Consider the automorphism $\varphi\in \Aut(F(X))$
defined by $\varphi(x) = x + [y,z]$, $\varphi $ acts identically
on the rest of variables from $X$. Then ${\hat f}_a(\varphi)(x)=
x+a[y,z]$. Therefore, the linearly stable automorphism ${\hat f}_a$
for $a\neq 1$ does not act identically on the group $\Aut(F(X))$.

2. Let  $X=\{x,y\}$. Then ${\Aut}_l(F(X))=\Aut(F(X))$ \cite{Co}.
In this case the automorphism ${\hat f}_a$ acts identically on the
group $\Aut(F(X))$. However, if $a\neq 1$ then ${\hat f}_a$ does
not act identically on the semigroup $\End(F(X))$. Indeed, take an
endomorphism $\varphi$ such that $\varphi(x)=[x,y]$. Then  ${\hat
f}_a(\varphi)(x)= a[x,y]$.
\end{remark}

\begin{proposition}
Let $|X|>2$. If $\xi$ is  a stable automorphism (acts identically
on $\Aut(F(X))$) then $\xi$ acts identically on $\End(F(X))$.
\end{proposition}

\begin{proof}
It follows from Remark \ref{R} that if $|X|>2$ then $\xi=\hat{f}_a$
is a stable automorphism if and only if $a=1$.
This proves the proposition.
\end{proof}

We  say that an automorphism $f$ of $\Aut(F(X))$ is an extendable
automorphism if there exists an automorphism $g$ of $\End(F(X))$
whose restriction to $\Aut(F(X))$ is $f$. It is obvious that all
extendable automorphisms of $\Aut(F(X))$ form a subgroup of
$\Aut(\Aut(F(X)))$. This subgroup is the image $Im(\tau)$ of a
homomorphism $\tau$ defined in the introduction.

\begin{corollary}
If $|X|>2$ then the group $\Aut(\End(F(X)))$ is isomorphic to the
subgroup of $\Aut(\Aut(F(X)))$ consisting of all extendable
automorphisms.
\end{corollary}

\section{Quasi-stable automorphisms of $\End(F(X))$}

In this section we define the notion of the quasi-stable automorphism
of $\End(F(X))$ and prove that any quasi-stable automorphism of
$\End(F(X))$ is inner.

Define, first, the diagonal automorphisms of the group ${\Aut}_l(F(X))$.
We proceed from the canonical isomorphism
$\delta:\Aut_l (F(X))\to\GL_n(P)$, $n=|X|$. Consider the commutative
diagram
$$
\CD
{\Aut}_l(F(X)) @>\delta>> \GL_n(P)\\
@. @/SE/ h  // @VV h_1 V\\
@. P^* \\
\endCD
$$
\noindent Define $\tilde {h}_1(A)=h_1(A)A,$ for every matrix $A\in
GL_n(P)$. Then $\tilde {h}_1$ is an automorphism of $GL_n(P)$. It
corresponds to the automorphism $\tilde {h}$ of  $\Aut_l (F(X))$
defined by $\tilde {h}= \delta^{-1}\tilde {h_1}\delta$. It is easy
to see that $\tilde {h}(g)(x)=h(g)g(x)$
 for every
$g\in\Aut_l (F(X))$ and every $x\in X.$ The automorphisms $\tilde {h}$ and
$\tilde {h_1}$ are called diagonal automorphisms.

\begin{definition}
An automorphism $\xi\in \Aut(\End(F(X)))$ is called quasi-stable if
the group ${\Aut}_l(F(X))$ is invariant in respect to $\xi$ and
 the restriction of $\xi$ to  ${\Aut}_l(F(X))$ is a diagonal
automorphism.
\end{definition}
Denote the restriction of  $\xi$ to  ${\Aut}_l(F(X))$ by
$\tau_0(\xi)$. Then $\tau_0(\xi)= \tilde {h}$ for some
homomorphism $h:\Aut_l (F(X))\to P^* $.

If $h(g)=1$ for every $g$ then the automorphism $\xi$ is linearly stable.
We show that every quasi-stable automorphism $\xi$ turns
out to be a linearly stable automorphism.

Consider the case  $X=\{x,y\}$. In this case  $\Aut
(F(X))={\Aut}_l(F(X))$ (see \cite{Co}). Thus, we do not need the
assumption that $\xi$ leaves ${\Aut}_l(F(X))$ invariant, and
$\tau_0(\xi)=\tau(\xi)$ where $\tau: \Aut(\End(F))\to
\Aut(\Aut(F))$ is the homomorphism defined in the introduction.

\begin{remark}
For the case $|X|>2$ we also could proceed from the homomorphism
$h:
\Aut (F(X))\to P^*$. However, in this case the corresponding
$\tilde h$ is not an automorphism of $\Aut (F(X))$.
\end{remark}

Let $\tau_0(\xi)$ be equal to $\tilde h$.
Denote the automorphism defined by the mapping $g(x)=y$, $g(y)=x$ and
$g(z)=z$ for the rest of elements of $X$ by $g_{xy}$. Denote
$h(g_{xy})=a_{xy}$. Notice that $g_{xy}^2=e$. Therefore,
$a_{xy}^2=1$. Since $P$ does not contain zero divisors $a_{xy}=\pm
1$. Denote $\varphi\in End_l(F(X))$ such that $\varphi(x)=x$ and
$\varphi(y)=a_{xy}x$ for any other $y\in X$ by $l$.

\begin{lemma}\label{Ch_F}
Let $\xi$ be a quasi-stable automorphism of $\End(F(X))$. Then
$\xi(c_u)=c_v l$, where $l\in End_l(F(X)$ is defined above. Thus,
$\xi(C_F)=C_F l$.
\end{lemma}

\begin{proof}
Let $g_{xy}\in \Aut(F)$ and $\xi(g_{xy})=a_{xy}g_{xy}$.
$c_u=c_ug_{xy}$. $\xi(c_u)=\xi(c_u)a_{xy}g_{xy}$. Let
$\xi(c_u)(x)=v$. Therefore, for any $y\in X$
$\xi(c_u)(y)=a_{xy}\xi(c_u)(x)=c_v l(y)$.
Thus, $\xi(c_u)=c_v l$.
\end{proof}

\begin{lemma}\label{Ch-l}
Let $\xi$ be a quasi-stable automorphism of $\End(F(X))$. Then for
any $x,y\in X$ $\xi(c_x)=c_x$ and $\xi(g_{xy})=g_{xy}$.
\end{lemma}

\begin{proof}
Let $c_x\in C_X$. It follows from the previous Lemma that
$\xi(c_x)=c_p l$ for some $p\in F(X)$. Let
$p=ax+b_1y_1+...+b_ky_k+p_1$, where $p_1$ is a non-linear
polynomial. $c_x$ is a right identity in the semigroup $C_F$ for
any $c_x\in C_X$. Therefore, $\xi(c_x)$ is a right identity for
$C_F$. Hence, $c_x\xi(c_x)=c_x$. Consequently,
$c_{(a+b_1+...+b_k)x}l=c_x$. Therefore, $a+b_1+...+b_k=1$ and $l$
is the identical mapping. Thus, $\xi(g_{xy})=g_{xy}$ for any
$x,y\in X$.

Let $g_{my}$ be the automorphism defined by the mapping
$g_{my}(x)=x$ and $g_{my}(y)=my$ for $y\neq x$, $y\in X$ and let
$h(g_{my})=d_m\neq 0$. Then $g_{my}c_x=c_x$. Therefore,
$\xi(g_{my}c_x)= d_mg_{my}c_p=\xi(c_x)=c_p$. Hence,
$d_mg_{my}(p)=p$. Remind that we present $p$ in the basis of Hall
of linearly independent elements. Thus, we obtain a system of
equalities of corresponding coefficients. $d_ma=a,
d_mmb_i=b_i,...$. If $a\neq 0$ then $d_m=1$.
Choosing $m\neq\pm 1$ we obtain that all other coefficients of $p$
are equal to $0$. In particular, $a=1$. Thus, in this case
$\xi(c_x)=c_x$.

Suppose now that $a=0$. Then $b_1+...+b_k=1$ and
$p=b_1y_1+...+b_ky_k+p_1$. Let $g'_{my}$ be the automorphism
generated by the mapping $g'_{my}(x)=x$ and $g'_{my}(y)=my+x$ for
$y\neq x$, $y\in X$ and let $h(g'_{my})=d'_m\neq 0$. Then
$g'_{my}c_x=c_x$. Therefore,
$\xi(g'_{my}c_x)=d'_mg'_{my}c_p=\xi(c_x)=c_p$. Hence,
$d'_mg'_{my}(p)=p$. Comparing the coefficients of $x$ we obtain
that $d'_m(b_1+...+b_k)=d'_m=0$. But this is impossible. Thus,
this case is impossible. Hence, $\xi(c_x)=c_x$.
\end{proof}

\begin{lemma}\label{Ch-lin}
Let $\xi$ be a quasi-stable automorphism of $\End(F(X)),$
$\tau_0(\xi)=\tilde h$.
Let $\End_l(F(X))$ be the subsemigroup of all linear endomorphisms.
Then $\xi(\End_l(F(X)))=\End_l(F(X))$.
\end{lemma}

\begin{proof}
Let $p$ be a linear polynomial. We proved above that
$\xi(c_p)=c_q$, $p,q\in F(X)$. Let $g_p$ be the automorphism
defined by the mapping $g(x)=p$ and $g(y)=y$ for $y\neq x$, $y\in
X$. Then $c_q=\xi(c_p)=\xi(gc_x)=h(g)g\xi(c_x)=
c_{h(g)g(x)}=c_{h(g)p}$. Hence, $q=h(g)p$ is a linear polynomial
as well. Let $\varphi\in End_l(F(X))$.
$\xi(c_{\varphi(x)})=\xi(\varphi c_x)=
\xi(\varphi)\xi(c_x)=\xi(\varphi)c_x=c_{\xi(\varphi)(x)}$. Since
$\varphi(x)$ is a linear polynomial, $\xi(\varphi)(x)$ is a linear
polynomial as well. Thus, $\xi(\varphi)\in End_l(F(X))$.
\end{proof}

The following Proposition can be obtained as a corollary of the
previous Lemma and the description of automorphisms of the linear
semigroup obtained by L.Gluskin in \cite{Gl} but we  give an
independent proof.

\begin{proposition}\label{ls-si}
Any quasi-stable  automorphism $\xi$ of $\End(F(X))$ is linearly stable.
Thus, $\xi$ is inner.
\end{proposition}

\begin{proof}
To prove the Proposition is enough to prove that a quasi-stable
automorphism acts identically on $End_l(F(X))$ and then use
Corollary \ref{lss}.

Let $p$ be a linear polynomial. Let $a$ be the sum of all
coefficients of $p$. It follows from the proof of Lemma
\ref{Ch-lin} that $\xi(c_p)=c_{kp}$, where $k\in P^*$. It follows
from the proof of Lemma \ref{Ch-l} that $h(g_{my})=d_m=1$, where
$g_{my}$ is the automorphism generated by the mapping
$g_{my}(x)=x$ and $g_{my}(y)=my$ for $y\neq x$, $y\in X$. Thus,
$\xi(c_xc_p)=\xi(c_{ax})=\xi(c_xg_{ay}c_y)=c_xg_{ay}c_y=c_{ax}$.
On the other hand
$\xi(c_xc_p)=\xi(c_x)\xi(c_p)=c_xc_{kp}=c_{akx}$. If $a\neq 0$
then $k=1$, that is $\xi(c_p)=c_{p}$. If $a=0$ then consider
$g_{my}$ such that the sum $a'$ of all coefficients of $g_{my}(p)$
is not equal to $0$. Then as we have proved above
$g_{my}c_p=\xi(g_{my}c_p)$. Therefore,
$g_{my}c_p=\xi(g_{my}c_p)=g_{my}\xi(c_p)$. Since $g_{my}$ is an
automorphism of $F(X)$ then $\xi(c_p)=c_p$.

Thus, $\xi$ acts identically on $C_l$. Let $\varphi\in
End_l(F(X))$. $c_{\varphi(x)}=\xi(c_{\varphi(x)})= \xi(\varphi
c_x)=\xi(\varphi)\xi(c_x)=\xi(\varphi)c_x=c_{\xi(\varphi)(x)}$.
Thus, $\xi(\varphi)(x)=\varphi(x)$. Hence,
$\xi(\varphi)=\varphi$.
\end{proof}

\section{Automorphisms of $\End(F(X))$}

In this section we
prove that any automorphism of $\End(F(x,y))$ is semi-inner.


\begin{theorem}\label{semi-inner}
Any automorphism of $\End(F(x,y))$ is a semi-inner automorphism.
\end{theorem}

\begin{proof}
P.Cohn \cite{C} proved that the group $\Aut(F(X))$ is generated by
linear and triangular automorphisms. Triangular automorphisms
assign $ax_i+f(x_1,...,x_{i-1})$ to $x_i$. Hence, for $X=\{x,y\}$
($X$ consists of two elements), a triangular automorphism assigns
$ay+f(x)$ to $y$. A Lie polynomial of one variable is a linear
polynomial. Therefore, the group $\Aut(F(x,y))$ consists of linear
automorphisms only. Thus, $\Aut(F(x,y))$ is isomorphic to
$GL_2(P)$. Let $\delta: \Aut(F(x,y))\rightarrow GL_2(P)$ be an
isomorphism. Then $\nu\rightarrow\delta^{-1}\nu\delta$ defines an
isomorphism $\Aut(GL_2(P))\rightarrow\Aut(\Aut(F(x,y)))$. If
$\nu\in\Aut(GL_2(P))$ is semi-inner then $\delta^{-1}\nu\delta$ is
a semi-inner automorphism of $\Aut(F(x,y))$. If
$\nu\in\Aut(GL_2(P))$ is diagonal then $\delta^{-1}\nu\delta$ is a
diagonal automorphism of $\Aut(F(x,y))$. It is well known
(\cite{OMe}) that the group of automorphisms of $GL_2(P)$ is
generated by semi-inner and diagonal automorphisms. Hence,  the
group of automorphisms of $\Aut(F(x,y))$ is generated by
semi-inner and diagonal automorphisms.

Let $\xi$ be an automorphism of $\End(F(X))$. In the introduction
we defined a homomorphism
$\tau:\Aut(\End(F(X)))\rightarrow\Aut(\Aut(F(X)))$, where
$\tau(\xi)=\xi^{\tau}$ is the restriction of $\xi$ to
$\Aut(F(X))$. Hence, $\tau(\xi)$ is a product of semi-inner and
diagonal automorphisms of $\Aut(F(X))$. Since a diagonal
automorphism commutes with a semi-inner automorphism we obtain
$\tau(\xi)=\xi_s\xi_d$, where $\xi_s$ is a semi-inner automorphism
and $\xi_d$ is a diagonal automorphism of $\Aut(F(X))$. $\xi_s$ is
a semi-inner automorphism of $\Aut(F(X))$ defined by a
semi-automorphism $f$ of $F(X)$. $f$ defines an automorphism
$\xi_1$ of $\End(F(X))$. Thus, $\xi_s=\xi_1^{\tau}$. Then
$\xi_1^{-1}\xi=\xi_2$ is an automorphism of $\End(F(X))$ and
$\xi_2^{\tau}=\xi_d$. Thus, $\xi_2$ is a quasi-stable
automorphism. It follows from Proposition \ref{ls-si} that $\xi_2$
is inner. Consequently $\xi=\xi_1\xi_2$ is semi-inner as a product
of semi-inner automorphisms.
\end{proof}

Now we formulate an application of Theorem 2. First, we consider
two problems. Let  $F_1$ and $F_2$ be free Lie algebras over a
field $P.$ Suppose that the semigroups of endomorphisms
$\End(F_1)$ and $\End(F_2)$ are isomorphic. Does this imply the
isomorphism of algebras $F_1$ and $F_2?$ Let $ \xi$ be an
isomorphism of $\End(F_1)$ and $\End(F_2)$. In which cases one can
state that there exists an isomorphism or a semi-isomorphism $f:
F_1\to F_2$ which induces $\xi, $ i.e., $\xi(\varphi)= f\varphi
f^{-}$ for every $\varphi\in \End(F_1)?$ If $f$ induces $\xi$ then
the pair $(f,\xi)$ defines the isomorphism or semi-isomorphism of
the actions of the semigroup $\End(F_1))$ on $F_1$ and
$\End(F_2))$ on $F_2$.

\begin{proposition}\label{iso}
If $F_1=F(x,y)$ then for any isomorphism $\xi: \End(F_1)\to\End(F_2)$
there exists a semi-isomorphism $f: F_1\to F_2$ which induces $\xi$.
\end{proposition}

\begin{proof}
Show that there exists an isomorphism $F_1\to F_2.$ We have to show
that if $F_2$ is freely generated by a set $Y$ then $ |Y|=2.$
Isomorphism $\xi$ induces an automorphism
of groups $\Aut(F_1)$ and $\Aut(F_2)$. The group
$\Aut(F_1)$ contains a non-trivial center consisting of scalar
automorphisms. If $|Y|\neq 2$ then the group $\Aut(F_2)$ does not
possess such a property. Thus, $|Y|=2$ and there is an isomorphism
$f_1: F_1\to F_2.$ We change this $f_1$ in order to get a desired
semi-isomorphism $f$. Define $\hat {f_1}: \End(F_1)\to \End(F_2)$ by the rule
$\hat {f_1}(\varphi)=f_1\varphi f_1^{-1}$ for every $\varphi \in \End(F_1)$.
The product $\hat {f_1}^{-1}\xi$ is an automorphism of the semigroup
$\End(F_1)$. Using Theorem 2 we get that this automorphism is semi-inner.
Thus, $\hat {f_1}^{-1}\xi=\hat g$, where $g$ is a semi-automorphism
of the algebra $F_1$. Now $\xi=\hat {f_1}\hat g$. Semi-isomorphism
$f=f_1g: F_1\to F_2$ induces the initial $\xi$.
\end{proof}

\begin{problem}
Does the theorem \ref{semi-inner} admit a generalization
for the case of arbitrary $X$, $|X|\geq 2$?
\end{problem}

\begin{proposition}\label{eq}
The following conditions on a free Lie algebra $F(X)$ are equivalent
\begin{enumerate}
\item
Any automorphism of $\End(F(X))$ is semi-inner.
\item
For any automorphism $\xi$ of $\End(F(X))$ the group
$\xi(Aut_l(F(X)))$ is conjugated to $Aut_l(F(X))$ (in the group
$\Aut(F(X))$).
\end{enumerate}
\end{proposition}

\begin{proof}
\begin{itemize}
\item
$1\Rightarrow 2$ There exists a semi-inner automorphism
$(\sigma,g)$ of $F(X)$  such that for any $\varphi\in \End(F(X))$
$\xi(\varphi)=g\sigma\varphi\sigma^{-1}g^{-1}$. For any $\alpha\in
Aut_l(F(X))$ we have $\sigma\alpha\sigma^{-1}\in Aut_l(F(X))$.
Therefore, $Aut_l(F(X))$ and $\xi(Aut_l(F(X)))$ are conjugated by
 $g\in\Aut(F(X))$.
\item
$2\Rightarrow 1$ Let $\xi$ be an automorphism of $\End(F(X))$.
$Aut_l(F(X))$ and $\xi(Aut_l(F(X)))$ are conjugated by
$g\in\Aut(F(X))$. $g$ defines an inner automorphism $\hat{g}$ of
$\End(F(X))$. $\xi\hat{g}^{-1}=\xi_1$ is an automorphism of
$\End(F(X))$ which induces an automorphism $\xi_2$ of
$Aut_l(F(X))$.  $\xi_2$ is semi-inner and, therefore, it is
extended to semi-inner automorphism $\hat{\xi}_2$ of $\End(F(X))$.
Automorphism $\xi_1\hat{\xi}_2^{-1}=\xi_3$ is a linearly stable
automorphism of $\End(F(X))$. Therefore, it is inner (see
Corollary \ref{lss}). Hence, $\xi=\xi_1\hat{g}=
\xi_3\hat{\xi}_2\hat{g}$ is semi-inner.
\end{itemize}
\end{proof}

\section{Reduction theorem}\label{6}

In this section we prove the general reduction theorem  for a
large class of varieties $\Theta$. This theorem  allows to reduce
the problem of the description of automorphisms of the category
$\Theta^0$ to the same problem  for a much simpler category
(consisting of two objects).

We assume that the variety $\Theta$  satisfies the following 3
conditions.

1. The variety $\Theta$ is hopfian. This means that every object
$F=F(X)$ of the category $\Theta^0$ is hopfian, i.e., every
surjective endomorphism $\nu: F\to F$ is an automorphism.

2. If $X=\{x_0\}$ is  a one element set and $F_0=F(x_0)$ is the
cyclic free algebra then for every automorphism $\varphi$ of the
category $\Theta^0$ we require $\varphi(F_0)$ is also a cyclic
free algebra $F(y_0)$.

3. We assume that there exists a finitely generated free algebra
$F^0=F(X^0)$, $X^0=\{x_1, \cdots, x_k\}$, generating
the whole variety $\Theta$, i.e., $\Theta=\Var \ (F^0)$.

For the sake of convenience in this paper we will
call a variety, satisfying these conditions, a hereditary variety.

We fix $F^0$ and $F_0$.

\begin{proposition} \label{cond}\cite{BPP}, (see also Appendix).
The conditions 1
and 2 imply that for every $F=F(X)$ and every $\varphi:
\Theta^0\to \Theta^0$ the algebras $F$ and $\varphi(F)$ are
isomorphic.
\end{proposition}

\begin{lemma} \cite{BPP}
 Any automorphism $\varphi: \Theta^0\to \Theta^0$ such that algebras $F$
and $\varphi(F)$ are isomorphic has the form
$$
\varphi=\varphi_0\varphi_1,
$$
\noindent
 where $\varphi_0$ is an inner automorphism of $\Theta^0$ and
$\varphi_1$ does not change objects.
\end{lemma}

Consider a constant morphism $\nu_0: F^0\to F_0$ such that
 $\nu_0(x)=x_0$ for every
$x\in X^0$.


\begin{theorem} ({\bf Reduction Theorem})\label{Red}
Let $\varphi$ be an automorphism of the category $\Theta^0$ which
does not change objects, and let $\varphi$ induce the identity
automorphism on the semigroup $\End(F^0)$ and
$\varphi(\nu_0)=\nu_0.$ Then $\varphi $ is an inner automorphism.
\end{theorem}

The proof of the Theorem consists of several steps.

1. It will be convenient to attach to the category $\Theta^0$ the
category of affine spaces $K^0_\Theta (H)$ over the algebra
$H=F^0$ \cite{Pl5}. The objects of $K^0_\Theta (H)$ have the form $\Hom
((F(X),H)$, where $F$ is an object of the category $\Theta^0$.
Morphisms
$$
\tilde s: \Hom ((F(X),H) \to \Hom ((F(Y),H)
$$
are defined by morphisms $s:F(Y)\to F(X)$ by the rule ${\tilde s}
(\nu)= \nu s$ for every $\nu: F(X)\to H$. We have a contravariant
functor $\Phi: \Theta^0\to K^0_\Theta (H)$.  The condition $\Var
(H)=\Theta$ implies that this functor yields the duality of the
categories $\Theta^0$ and $K^0_\Theta(H)$ (see \cite{Pl5} and
Appendix). Consider the automorphism $\varphi^H$ of the category of
affine spaces which is the image of $\varphi$ under the duality
above.  Functor $\varphi^H: K^0_\Theta(H)\to K^0_\Theta(H)$ does not
change objects and for $s: F(Y)\to F(X)$ we define $\varphi^H(\tilde
s)=\widetilde{\varphi(s)}$. This definition is correct, since
$\tilde{s_1}=\tilde{s_2}$ implies $s_1=s_2$.

We will show that $\varphi^H$ is in a certain sense 
a quasi-inner automorphism. First of all, $\varphi$ defines a
substitution on each set $\Hom(F(X),H)$. Indeed,  $\nu: F(X)\to H$
and $\varphi(\nu): F(X)\to H$ give rise to a substitution $\mu_X$
defined by  $\mu_X(\nu)=\varphi(\nu).$ The following proposition
explains the transition from $\Theta^0$ to $K^0_\Theta(H)$.

\begin{proposition} Let $s:F(Y)\to F(X)$. Then
$$
\varphi^H(\tilde {s})=\mu_Y{\tilde s}\mu_X^{-1}: \Hom(F(X),H)\to
\Hom(F(Y),H).
$$
\end{proposition}

\begin{proof}
For every $s: F(Y)\to F(X)$ and every  $\nu: F(X)\to H$ the
equality $\varphi^H(\tilde s)(\nu)=\widetilde{\varphi(s)}(\nu)=
\nu \varphi(s)$ holds. Therefore, we have
$$
\mu_Y\tilde s\mu_X^{-1}(\nu)=\mu_Y(\mu_X^{-1}(\nu)s)=
\varphi(\varphi^{-1}(\nu)s)=\nu\varphi(s)=\varphi^H(\tilde s)(\nu).
$$
\end{proof}

\begin{remark} If the automorphism $\varphi^H$ has a presentation above
we call it quasi-inner.
\end{remark}

Consider separately the case $X=X^0$ and take the substitution
$\mu_{X^0}:\Hom(F^0,H)\to\Hom(F^0,H)$. By the condition of the
theorem the equality $\mu_{X^0}(\nu)=\varphi(\nu)=\nu$ holds for
any $\nu: F^0\to H=F^0$. This means that $\mu_{X^0}=1$. Then for
$s: F(Y)\to F(X^0)$ we have
$$\varphi^H(\tilde s)=\mu_Y\tilde s\mu_{X^0}^{-1}=\mu_Y\tilde s
= \widetilde{\varphi(s)}.
$$

For $s: F(X^0)\to F(Y)$ we get
$$
\varphi^H(\tilde s)=\mu_{X^0}\tilde s\mu_{Y}^{-1}=\widetilde
s\mu_Y^{-1} = \widetilde{\varphi(s)}.
$$
Therefore, $\tilde s=\widetilde{\varphi(s)}\mu_Y$.

2. Now we use the category of polynomial maps $\Pol_\Theta(H)$.
Objects of this category have the form $H^n$, where $n$ changes
and
$H$ is fixed. Morphisms are represented by polynomial maps
$s^\alpha : H^n\to H^m$ defined below. Take a set
$X=\{x_1,\cdots,x_n\}$. Denote $\alpha_X:\Hom(F(X),H)\to H^n$ the
canonical bijection defined by
$\alpha_X(\nu)=(\nu(x_1),\cdots,\nu(x_n))$ for every $\nu:F(X)\to
H$. Let now $s:F(Y)\to F(X)$ be given and $X=\{x_1,\cdots,x_n\}$
$Y=\{y_1,\cdots,y_m\}$. Consider the diagram

$$\CD \Hom(F(X),H) @>\tilde s>> \Hom(F(Y),H) \\ @V\alpha_X VV
@VV\alpha_Y V\\
H^n@>s^\alpha>> H^m\endCD$$ \noindent where
$s^\alpha=\alpha_Y\tilde s\alpha_X^{-1};$ $\tilde s=
\alpha^{-1}_Ys^\alpha\alpha_X$. Then $s^\alpha(a_1,\cdots,a_n)=
\alpha_Y\tilde s\alpha_X^{-1}(a_1,\cdots,a_n)$. Take a point
$\nu=\alpha^{-1}_X(a_1,\cdots,a_n): F(X)\to H$. Then
$$
s^\alpha(a_1,\cdots,a_n)=\alpha_Y\tilde s(\nu)=\alpha_Y(\nu s)=
(\nu s(y_1),\cdots,\nu s(y_m)).
$$

Denote $s(y_i)=w_i(x_1,\cdots,x_n)$, $i=1,\cdots,m$. We have got
$$
s^\alpha(a_1,\cdots,a_n)=(w_1(a_1,\cdots,a_n),\cdots,
w_m(a_1,\cdots,a_n)).
$$
Indeed,
$$
s^\alpha(a_1,\cdots,a_n)=(\nu(w_1(x_1,\cdots,x_n)),\cdots,
\nu(w_m(x_1,\cdots,x_n)))=
$$
$$
=(w_1(x^\nu_1,\cdots,x^\nu_n),\cdots,
w_m(x^\nu_1,\cdots,x^\nu_n))= (w_1(a_1,\cdots,a_n),\cdots,
w_m(a_1,\cdots,a_n)).
$$
Thus, we defined morphisms $s^\alpha:  H^n\to H^m$ in the category
$\Pol_\Theta(H)$.

Consider special constant morphisms in the category $\Theta^0$.
First, take  morphisms of the form: $\nu=\nu_a: F_0\to F(X)$
defined by $\nu_a(x_0)=a$, $a\in F(X)$. Recall that the constant
morphism $\nu_0: F^0\to F(x_0)$ is defined by $\nu_0(x)=x_0$ for
every $x\in X^0$.

Take $\nu=\nu_a\nu_0: F(X^0)\to F(X)$. Then
 $\nu(x)=a$ for every $x\in X^0$, and $\nu$ is a constant we will
be dealing with.

Let, further, $\varphi$ be an automorphism of $\Theta^0$ which
does not change objects. This $\varphi$ induces a substitution on
each set $\Hom(F(X),F(Y))$ denoted by $\mu_{X,Y}$. In particular,
$\mu_{X,X^0}=\mu_X$. The substitution $\mu_{x_0,X}$ on the set
$\Hom(F_0,F(X))$ induces the substitution $\sigma_X$ on the
algebra $F(X)$ defined by the rule
$\varphi(\nu_a)=\nu_{\sigma_X(a)}$. It is proved \cite{Pl5} that
for every $\mu: F(X)\to F(Y)$ the formula
$$
\varphi(\mu)=\sigma_Y\mu\sigma_X^{-1}
$$
holds. In this sense the automorphism $\varphi$ is said to be
a quasi-inner automorphism in the category $\Theta^0$. Take now
$\nu=\nu_a\nu_0$. Then
$$\varphi(\nu)=\varphi(\nu_a)\varphi(\nu_0)=\nu_{\sigma_X(a)}\varphi(
\nu_0).$$ If $\varphi$ does not change $\nu_0$ then
$$\varphi(\nu)=\nu_{\sigma_X(a)}\nu_0.$$
For every $x\in X^O$ we get $\varphi(\nu)(x)=\sigma_X(a)$, here
$\varphi(\nu)$ is also a constant.

Now we are in the position to make the next step. We return to the
category of polynomial maps and consider how the constant maps
defined above in $\Theta^0$ look like in $\Pol_\Theta(H)$. Take
$s: H=F(X^0)=F^0\to F(X)$ defined by $s=s_w\nu_0$, $w\in F(X)$.
Then $s(x)=w$ for every $x\in X^0$.  Let $X^0=\{x_1,\cdots,
x_k\}$. We get the commutative diagram

$$\CD \Hom(F(X),H) @>\tilde s>> \Hom(F(X^0),H) \\ @V\alpha_X VV
@VV\alpha_{X^0} V\\
H^n@>s^\alpha>> H^k\endCD$$ \noindent

Then $ s^\alpha(a_1,\cdots,a_n)=(w(a_1,\cdots,a_n),\cdots,
w(a_1,\cdots,a_n)) $, where $w(a_1,\cdots,a_n)$ is taken $k$
times. Considering the projection $\pi: H^k\to H$,
$\pi(b_1,\cdots,b_k)=b_1$ we get $\pi s^\alpha(a_1,\cdots,a_n)=
w(a_1,\cdots,a_n)$.

Take an arbitrary $s:F(Y)\to F(X)$. Let $X=\{x_1,\cdots,x_n\}$,
and $Y=\{y_1,\cdots,y_m\}$. Let $s(y_i)=w_i(x_1,\cdots,x_n)=w_i$.
Take a constant map $s_i=\nu_{w_i}\nu_0: F(X^0)\to F(X)$. The
sequence $s_1,\cdots, s_m$ depends on $s $ and on the basis of
$Y$. In this situation we denote
$$
    s=_Y(s_1,\cdots, s_m).
$$
We have also $ s^\alpha: H^n\to H^m$, and $s_i^\alpha: H^n\to
H^k$. There is a relation between $s^\alpha$ and $s_i^\alpha$,
$i=1,\cdots,m$:
$$
 s^\alpha(a_1,\cdots,a_n)= (\pi
s_1^\alpha(a_1,\cdots,a_n),\cdots, \pi
s_m^\alpha(a_1,\cdots,a_n)).$$
 Indeed,
$$
s^\alpha(a_1,\cdots,a_n)= (w_1(a_1,\cdots,a_n), \cdots,
w_m(a_1,\cdots,a_n))= (\pi s_1^\alpha(a_1,\cdots,a_n), \cdots ,
\pi s_m^\alpha(a_1,\cdots,a_n)).
$$

This formula is a key working tool for the proof of the theorem.
It was the reason to replace the category of affine spaces by the
category of polynomial maps.

Now we are able to prove the reduction theorem.

{\sc Proof of Theorem \ref{Red}.}
Let us return to the automorphism $\varphi: \Theta^0\to \Theta^0$.
For every algebra $F=F(X)$, $X=\{x_1,\cdots,x_n\}$ we will
construct an automorphism $\sigma_X: F\to F$ depending on
$\varphi$. The collection of such automorphisms will define
$\varphi$ as an inner automorphism.

Consider morphisms $\varepsilon_i=\nu_{x_i}\nu_0: F^0\to F$,
$i=1,\cdots,n$. We have
$\varphi(\varepsilon_i)=\varphi(\nu_{x_i})\nu_0$, and let
$\varphi(\nu_{x_i})(x_0)=y_i=\varphi(\varepsilon_i)(x)$ for every
$x\in X^0$. Let $Y=\{y_1,\cdots,y_n\}$. From the proof of
Proposition \ref{cond} follows that if the variety $\Theta$ is hopfian
then $Y$ is also a basis in $F$.

 Define the automorphism
$\sigma_X: F\to F$ by the rule $\sigma_X(x_i)=y_i$.

 Let $s$ be an
automorphism of the algebra $F=F(X)$, and let $s(x_i)=
w_i(x_1,\cdots,x_n)
=w_i$. Take $\nu_{w_i}: F_0\to F$, $\nu_{w_i}=
s\nu_{x_i}$ and let
$s_i=\nu_{w_i}\nu_0=s\nu_{x_i}\nu_0=s\varepsilon_i$.

We have $s=_X(s_1,\cdots,s_n)$. We will check that
$\varphi(s)=_Y(\varphi(s_1),\cdots,\varphi(s_n)$.

We have to verify that if $\varphi(s)(y_i)=w_i'$ then
$\varphi(s_i)(x)=w_i'$ for every $x\in X^0$. Compute
$$
\varphi(s)(y_i)=\varphi(s)\varphi(\nu_{x_i})(x_0)=
\varphi(s\nu_{x_i})(x_0)= \varphi(\nu_{w_i})(x_0)=
$$
$$
= \varphi(\nu_{w_i})\nu_0(x)=\varphi(\nu_{w_i})\varphi(\nu_0)(x)=
\varphi(\nu_{w_i}\nu_0)(x)=\varphi(s_i)(x)
$$
for every $x\in X^0$. Thus, $\varphi(s)(y_i)=w_i'=\varphi(s_i)(x)$
for every $x\in X^0$.

This implies $\varphi(s)\sigma_X=_X (\varphi(s_1),\cdots,
\varphi(s_n)$. Indeed, $\varphi(s)\sigma_X(x_i)=\varphi(s)(y_i)=
\varphi(s_i)(x)$ for every $x\in X^0$.



Consider the image of the formula $\tilde
{s_i}={\widetilde{\varphi(s_i)}}\mu_X$, $i=1,\cdots,n$ in the
category of polynomial maps $\Pol_\Theta(H)$.

Take the diagram
$$\CD \Hom(F(X),H) @>\mu_X>> \Hom(F(X),H) \\ @V\alpha_X VV @VV\alpha_{X}
V\\
H^n@>\mu_X^\alpha>> H^n\endCD$$

We have got a map $\mu_X^\alpha=\alpha_X\mu_X\alpha_X^{-1}: H^n\to
H^n$. In particular,
$\mu_{X^0}^\alpha=\alpha_{X^0}\mu_X\alpha_{X^0}^{-1}$. By the
condition of  Theorem, $\mu_{X^0}=1$ and $\mu^\alpha_{X^0}=1$.

Since, $\varphi^H(\tilde s)=\mu_X\tilde s \mu_X^{-1}=
\widetilde{\varphi(s)}$ for $\varphi(s):F(X)\to F(X)$ the following equality
holds
$$
\varphi(s)^\alpha=\alpha_X\widetilde{\varphi(s)}\alpha^{-1}_X=
\alpha_X\mu_X\alpha_X^{-1}\alpha_X\tilde s
\alpha_X^{-1}\alpha_X\mu_X^{-1}\alpha^{-1}_X=\mu^\alpha_X
s^\alpha(\mu^\alpha_X)^{-1}.
$$
For $s_i: F^0\to F$ we get $\tilde
{s_i}=\widetilde{\varphi(s_i)}\mu_X$ and hence,
$$
s^\alpha_i=\alpha_{X^0}\tilde {s_i} \alpha_X^{-1}=
\alpha_{X^0}\widetilde{\varphi(s_i)}\alpha_X^{-1}\alpha_X\mu_X\alpha_X^{-1}=
\varphi(s_i)^\alpha \mu_X^\alpha,
$$
where $s^\alpha_i=\varphi(s_i)^\alpha \mu_X^\alpha$ are polynomial
mappings from  $H^n$ to $H^k.$ For $a=(a_1,\cdots,a_n)\in H^n$ we
have
$$
s^\alpha(a_1,\cdots,a_n)= (\pi
s_1^\alpha(a_1,\cdots,a_n),\cdots,\pi s_n^\alpha(a_1,\cdots,a_n))=
$$
$$
=(\pi\varphi(s_1)^\alpha\mu_X^\alpha(a_1,\cdots, a_n),\cdots,
\pi\varphi(s_n)^\alpha\mu_X^\alpha(a_1,\cdots, a_n)).
$$
Take
$$\varphi(s)\sigma_X=_X(\varphi(s_1),\cdots,\varphi(s_n)$$
 and apply this formula to the point $\mu^\alpha_X(a_1,\cdots,a_n)$.
Then
$$
(\varphi(s)\sigma_X)^\alpha(\mu_X^\alpha(a_1,\cdots,a_n))=
(\pi\varphi(s_1)^\alpha\mu_X^\alpha(a_1,\cdots, a_n),\cdots,
\pi\varphi(s_n)^\alpha\mu_X^\alpha(a_1,\cdots, a_n))=s^\alpha(a).
$$

 Thus, $s^\alpha(a)=(\varphi(s)\sigma_X)^\alpha\mu_X^\alpha(a)$.
 Since this formula holds in every point $a$ then
 $$s^\alpha=(\varphi(s)\sigma_X)^\alpha\mu_X^\alpha=
 \sigma_X^\alpha\varphi(s)^\alpha\mu_X^\alpha.$$
 Hence,
 $\mu_X^\alpha=(\varphi(s)^{-1})^\alpha(\sigma_X^{-1})^\alpha
 s^\alpha= (s\sigma_X^{-1}\varphi(s)^{-1})^\alpha$.

 Denote $\xi_X=s\sigma_X^{-1}\varphi(s)^{-1}$. This is an
 automorphism  of the algebra $F=F(X)$ and
 $\mu_X^\alpha=\xi_X^\alpha$. Therefore, $\tilde {\xi_X}=\mu_X$.
 In particular, $\xi_X$ does not depend on the choice of the
 automorphism $s$.

 Let now an arbitrary  $\delta: F(X)\to F(Y)$ be given. Then we have
 $$\varphi^H(\tilde \delta)=\mu_X\tilde \delta\mu_Y^{-1}=
 \tilde
{\xi_X}\tilde\delta{\tilde\xi}_Y^{-1}=\widetilde{\xi_Y^{-1}\delta\xi_X}=
 \widetilde{\varphi(\delta)}
 $$
 This gives $\varphi(\delta)=\xi^{-1}_Y\delta\xi_X$.

 Since our initial $s$ is
 arbitrary,
 one can take $s=1$. Then $\xi_X=\sigma^{-1}_X$.


 Finally we  get
 $$
 \varphi(\delta)=\sigma_Y\delta\sigma_X^{-1}
 $$
\hfill $\Box$


\section{Automorphisms of the category of free Lie algebras. The proof
of the main Theorem}

Return to the variety $\Theta$ of all Lie algebras over an
infinite field. We want to prove that every automorphism of the
category $\Theta^0$ is semi-inner.


{\sc Proof of Theorem 1.}
 This variety $\Theta$ is hopfian and  is generated by the free Lie algebra
$F^0=F(x,y)$ \cite{BK}. It is clear that condition 2 is also
valid. Thus, the conditions  from section \ref{6} are fulfilled.
Therefore, the variety $\Theta$ is hereditary.

 It is enough to consider automorphisms $\varphi$
which do not change objects \cite{MPP}. Take such $\varphi$ and
induce the automorphism $\varphi_{F^0}$ of the semigroup
$\End(F^0)$. According to Theorem \ref{semi-inner} such
automorphism is semi-inner and is defined by the semi-automorphism
$(\sigma,s_{F^0}): F^0\to F^0$. For every algebra $F=F(X)$, which
is distinct from $F^0$ take a semi-automorphism
$(\sigma,\sigma_F): F\to F$.
Semi-automorphisms $(\sigma,s)_{F^0}=(\sigma, s_{F_0})$ and
$(\sigma,s)_{F}=(\sigma,\sigma_F)$ define a semi-inner automorphism
$\psi$ of the category $\Theta^0$.
This $\psi$ does not change objects.
Automorphisms $\varphi$ and $\psi$ act in
the same way on the semigroup $\End(F^0)$. Thus, the automorphism
$\varphi_1=\psi^{-1}\varphi$ acts on $\End(F^0)$ identically.

Take a constant morphism $\nu_0: F^0\to F_0$ with
$\nu_0(x)=\nu_0(y)=x_0$. Let us verify that $\varphi_1(\nu_0)$ is also a
constant. Take an automorphism $\eta$ of the algebra $F^0$ defined
by $\eta(x)=y$, $\eta(y)=x$. We have $\nu_0\eta=\nu_0$. Therefore
$\varphi_1(\nu_0\eta)=\varphi_1(\nu_0)\eta= \varphi_1(\nu_0).$
Hence,
$\varphi_1(\nu_0)(x)=\varphi_1(\nu_0)\eta(x)=\varphi_1(\nu_0)(y)=ax_0$
for $a\neq 0$.

Automorphisms of free Lie algebras: ${f_{F_0}}(x_0)=ax_0$ and $f_F(x)=x$
for $x\in X$, $F=F(X)\neq F_0$, define an inner automorphism $\hat f$
of the category $\Theta^0$, which does not  change objects.
Observe that isomorphism $f_F$ acts trivially on $F^0$.
 We have $\hat{f}(\nu_0)=f_{F_0}\nu_0 f_{F^0}^{-1}$ and
$$
\hat{f}(\nu_0)(x)={f}_{F_0}\nu_0 {f}_{F^0}^{-1}(x)=
{f}_{F_0}\nu_0(x)={f}_{F_0}(x_0)=ax_0,
$$
 Thus, $\varphi_1(\nu_0)$ and $\hat{f}(\nu_0)$ coincide.
Therefore ${\hat{f}}^{-1}\varphi_1(\nu_0)=\nu_0$. Denote
${\hat{f}}^{-1}\psi^{-1}\varphi=\psi_0$. Then
$\psi_0(\nu_0)=\nu_0$ and $\psi_0$ acts trivially on $\End(F(x,y))$.
 By Reduction Theorem the automorphism $\psi_0$
is inner.  We have got
$\varphi=\psi\hat{f}\psi_0=\hat \sigma \psi_1\hat f\psi_0.$

Thus, $\varphi$ is semi-inner and the theorem is proved.
\hfill $\Box$

Along with the automorphisms of categories of free algebras of
varieties it is natural to consider also the autoequivalences of
these categories (see Section 1). Let $(\varphi,\psi) $ be an
autoequivalence of the category of free Lie algebras. We call it
semi-inner if the functors $\varphi$ and $\psi$ are
semi-isomorphic to the identity functor.

It was proved in \cite{Z},  that for every $\Theta$ and every
autoequivalence $(\varphi, \psi)$ of the category $\Theta^0$ there
are factorizations
$$
\varphi=\varphi_0\varphi_1, \qquad\qquad \psi=\varphi_1^{-1}\psi_0,
$$
\noindent where $\varphi_0$ and $\psi_0$ are isomorphic to the identity
functor and $\varphi_1$ is an automorphism.

This means that every autoequivalence is isomorphic to an
automorphism. This remark and Theorem 1 lead to the following
statement:

\begin{theorem}\label{Eq}
Every autoequivalence of the category of free Lie
algebras is semi-inner.
\end{theorem}

In the introduction we discussed the categories of algebraic sets
$K_\Theta(H)$, $H\in \Theta$. The two problems were pointed out,
namely, the problem of isomorphism of categories $K_\Theta(H_1)$ and
$K_\Theta(H_2)$ and the problem of equivalence of the same categories.
 For the variety of Lie algebras and algebras $H_1$ and $H_2$, satisfying
$\Var(H_1)=\Var(H_2)=\Theta$, the first problem is solved in \cite{Pl6}
with the aid of Theorem 1, while the solution of the second problem
also in \cite{Pl6} requires arguments from Theorem \ref{Eq}.

\section*{Appendix}

In this section we prove two propositions we have referred to in the
text.

Remind that contravariant functor $\Phi: \Theta^0\to K^0_\Theta (H)$
assign $\Hom ((F(X),H)$ to any $F(X)\in\Theta^0$ and for any
$s:F(Y)\to F(X)$ it assigns $\tilde s: \Hom ((F(X),H) \to \Hom ((F(Y),H)$
defined by the rule ${\tilde s}(\nu)= \nu s$.

\begin{proposition} The functor $\Phi: \Theta^0\to K^0_\Theta (H)$
defines a duality of categories if and only if $\Var(H)=\Theta$.
\end{proposition}

\begin{proof}
In our case the duality of categories means that if
$s_1$, $s_2$ morphisms $F(Y)\to F(X)$ then $s_1\neq s_2$ implies
$\widetilde s_1\neq \widetilde s_2$.
Let $s_1\neq s_2$ and assume $\widetilde s_1=\widetilde s_2$.
Take $y\in Y$ such that $s_1(y)=w_1$, $s_2(y)=w_2$,
and $w_1\neq w_2$. Show that the non-trivial identity $w_1\equiv w_2$
is fulfilled
in algebra $H$. Take an arbitrary $\nu: F(X)\to H$. The equality
$\widetilde s_1=\widetilde s_2$ yields $\widetilde s_1(\nu)=
\widetilde s_2(\nu)$.   We have also $\nu s_1=\nu s_2$. Apply this
equality to $y$. We get  $\nu s_1(y)=\nu(w_1)=\nu s_2(y)=\nu (w_2).$
Since $\nu$ is arbitrary, we get that $\widetilde s_1= \widetilde s_2$
implies  $w_1\equiv w_2$ in $H$.

Assume that  $\Var(H)=\Theta$. Then there are no non-trivial identities in
$H$. This means that the equality $\widetilde s_1= \widetilde s_2$
does not hold in $K^0_\Theta (H)$. We proved that if $\Var(H)=\Theta$
then $s_1\neq s_2$ implies  $\widetilde s_1\neq \widetilde s_2$
and we get a duality of categories.

Conversely, let us prove that if  $\Var(H) \neq\Theta$ then there
is no duality. Since   $\Var(H) \neq\Theta$ there exists a
non-trivial identity $w_1\equiv w_2$ in $H$, where $w_1, w_2$ in
some $ F(X).$ Take $Y=\{y_0\}$. Consider $s_1$ and $s_2$ from
$F(Y)$ to $F(X)$ defined by the rule: $s_1(y_0)=w_1$,
$s_2(y_0)=w_2$. Show that $\widetilde s_1= \widetilde s_2.$ This
will mean that there is no duality. Take an arbitrary $\nu:
F(X)\to H$. Then $\widetilde s_1(\nu)=\nu s_1$, $\widetilde
s_2(\nu)=\nu s_2$, both $F(Y)\to H$. Take $y_0$. Then $\nu
s_1(y_0)=\nu(w_1)$, $\nu s_1(y_0)=\nu(w_2)$. Since  $w_1\equiv
w_2$ is an identity in $H$ then $\nu(w_1)=\nu(w_2)$ and
correspondingly,  $\nu s_1(y_0)=\nu s_2(y_0)$. Since the set $Y$
consists of one element $y_0$ then $\nu s_1 =\nu s_2$. This
equality holds for every $\nu$ and therefore,  $\widetilde
s_1=\widetilde s_2$.
\end{proof}

\begin{proposition} \label{hopf} Let the variety $\Theta$ be hopfian,
$\varphi$ an automorphism of $\Theta^0$ and $\varphi(F_0)$ is
isomorphic to $F_0$. Then $\varphi(F)$ is isomorphic to $F$ for every $F=F(X)$.
\end{proposition}

We use some new notions to prove Proposition \ref{hopf}.

\begin{definition} Let $X$ be a set in a free algebra $F=F(Y)$. We say
that $X$ defines freely algebra $F$ if every map $\mu_0: X\to F$
can be extended uniquely up to endomorphism $\mu: F\to F$.
\end{definition}

\begin{remark} In many cases the notions ``to define freely'' and
``to generate freely'' coincide. For instance, this is true for the
variety of all groups (E. Rips, unpublished). If this is true for the variety
of Lie algebras we do not know.
\end{remark}

\begin{lemma}
Let the variety $\Theta$ be hopfian. Let $|X|\geq |Y|$. Then $X$
defines freely $F=F(Y)$ if and only if $X$ is a basis in $F$
and $|X|=|Y|$.
\end{lemma}

\begin{proof}
Take an arbitrary surjection $\mu_0: X\to Y$. If $X$ defines
freely $F$ then there exists surjective endomorphism $\mu: F\to
F$. Since $F$ is hopfian, $\mu$ is automorphism. Then $\mu_0$ is a
bijection. The inverse bijection defines the inverse automorphism.
Therefore,
 $|X|=|Y|$ and $X$ is a basis in $F$. The ``only if'' part is evident.
\end{proof}

For the sake of self-completeness we repeat some material from
\cite{BPP}. Take a free algebra $F=F(X)$ and consider a system
of morphisms $\epsilon_i: F_0\to F$, $i=1,\cdots, n$.

\begin{definition} A system of morphisms  $(\epsilon_1,\cdots, \epsilon_n)$
defines freely an algebra $F$ if for every sequence of homomorphisms
$f_1, \cdots, f_n$, $f_i: F_0\to F$ there exists a
unique endomorphism
$s: F\to F$ such that $f_i=s\epsilon_i$ where  $i=1,\cdots, n$.
\end{definition}

It is proved in \cite{BPP} that the system $(\epsilon_1,\cdots, \epsilon_n)$
defines freely $F$ if and only if the system of elements
$(\epsilon_1(x_0),\cdots, \epsilon_n(x_0))$ defines freely algebra $F$.
It is obvious, that if the system $(\epsilon_1,\cdots, \epsilon_n)$
defines freely $F$  then the system $(\varphi(\epsilon_1),\cdots,
\varphi(\epsilon_n))$ defines freely $\varphi(F)$ if $\varphi(F_0)=
F(y_0)$.

{\sc Proof of Proposition \ref{hopf}.}
Let the variety $\Theta$ be hopfian, $\varphi(F_0)=F(y_0)$. Let
$\varphi(F)=F(Y)$ where $F=F(X)$. We prove that algebras $F(X)$
and $F(Y)$ are isomorphic.

Let, first,
$|X|\geq |Y|$ and $X=\{ x_1, \cdots   , x_n\}$. Define the system
 ($\epsilon_1,\cdots, \epsilon_n)$ by the condition
 $\epsilon_i(x_0)=x_i$, for every $i$. This system defines freely
algebra $F$. Then the system
$(\varphi(\epsilon_1),\cdots,
\varphi(\epsilon_n))$ defines freely algebra $\varphi(F)=F(Y)$.
This means that the set $Y'$ of elements $y'_i=\varphi(\epsilon_i)(y_0)$
defines freely algebra $F(Y)$. Since $|Y'|=|X|\geq |Y|$ then the
system $Y'$ is a basis in $F(Y)$ and   $|Y'|=|X|= |Y|$. The map
$x_i\to y'_i$ defines the isomorphism of algebras $F(X)$ and $F(Y)$.

Let now $|X| < |Y|$. Then $F(X)=\varphi(F(Y))^{-1}$. Applying the
same method to the automorphism $\varphi^{-1}$ we get the
contradiction.
\hfill $\Box$


\end{document}

\bye